\documentclass[reqno,11pt]{amsart}

\usepackage{geometry}
\usepackage{mathtools,amssymb,amsthm}
\usepackage[T1]{fontenc}
\usepackage[utf8]{inputenc}
\usepackage{color,enumitem,needspace,etoolbox}
\usepackage[
  colorlinks=true,
  linkcolor=blue,
  anchorcolor=green,
  citecolor=blue
]{hyperref}

\definecolor{bleu1}{RGB}{0,57,128}
\patchcmd{\section}{\normalfont}{\normalfont\color{bleu1}}{}{}
\patchcmd{\subsection}{\normalfont}{\normalfont\color{bleu1}}{}{}
\patchcmd{\subsubsection}{\normalfont}{\normalfont\color{bleu1}}{}{}

\newtheorem{proposition}{Proposition}[section]
\newtheorem{theorem}{Theorem}[section]
\newtheorem{lemma}{Lemma}[section]
\newtheorem{corollary}{Corollary}[section]

\hypersetup{
  pdftitle={Sharpness of the acceleration-dependent Holder exponent for the integrated density of states},
  pdfsubject={Complete revised manuscript},
  pdfauthor={Jiabo Ye and Jiangong You}
}

\begin{document}

\title[Sharpness of IDS H\"older exponents]{Sharpness of the acceleration-dependent H\"older exponent for the integrated density of states}
\author{Jiabo Ye}
\author{Jiangong You}
\date{}

\begin{abstract}
We prove that the acceleration-dependent H\"older exponent $1/(2\omega)$ for the integrated density
of states is sharp in the supercritical regime. For every positive integer $k$ and every Diophantine
frequency, we construct a real even trigonometric polynomial potential of degree $k$ such that, at a
spectral edge, the Lyapunov exponent is positive, the acceleration is $k$, and the exact local
H\"older exponent is $1/(2k)$. In particular, the $1/2$-H\"older regularity of the subcritical
regime for Diophantine frequencies does not hold in general in the supercritical regime. The proof
uses Aubry duality and a KAM iteration that adjusts the potential coefficients to make the dual
cocycle reducible to a principal unipotent matrix. A comparison with an explicit reference family
and symplectic block diagonalization then reduce the IDS estimate to a two-dimensional
rotation-number estimate.
\end{abstract}

\maketitle

\section{Introduction}

\subsection{Background and motivation}

We consider the one-frequency Schr\"odinger operators
\begin{equation}\label{e:op}
  (H_{v,\theta}u)_{n}=u_{n+1}+u_{n-1}+v(\theta+n\alpha)u_{n}, \qquad n\in\mathbb{Z},
\end{equation}
on $\ell^{2}(\mathbb{Z})$, where $v:\mathbb{T}\to\mathbb{R}$ is analytic,
$\mathbb{T}=\mathbb{R}/\mathbb{Z}$, and $\alpha$ is irrational. The spectrum $\Sigma_{v}$ is
independent of $\theta$. The integrated density of states (IDS) is defined by
\begin{equation*}
  N_{v}(E)=\int_{\mathbb{T}} \langle\delta_{0},\mathbf{1}_{(-\infty,E]}(H_{v,\theta})\delta_{0}\rangle\,\mathrm{d}\theta
\end{equation*}
For background on ergodic Schr\"odinger operators and their
spectral theory, we refer to Damanik--Fillman \cite{DF22}.

For small analytic potentials and Diophantine frequencies, Hadj Amor \cite{HA06,HA09} proved
$1/2$-H\"older continuity of the IDS. Avila--Jitomirskaya \cite{AJ10} obtained the same regularity
with a smallness threshold independent of the frequency.

In the positive Lyapunov exponent regime, Bourgain \cite{B00}
proved H\"older estimates for the almost Mathieu operator at sufficiently
large coupling and Diophantine frequencies. Goldstein--Schlag \cite{GS}
established H\"older continuity for general analytic potentials under a
strong Diophantine condition, and subsequently obtained finer estimates
for analytic perturbations of trigonometric polynomial potentials
\cite{GS08}. You--Zhang \cite{YZ14} extended H\"older regularity to all
Diophantine and sufficiently weak Liouville frequencies.
Further refinements by Han--Zhang \cite{HZ22} and Han--Schlag
\cite{HSL26} relate the admissible rational approximation rate
to the Lyapunov exponent. In the large-coupling regime, Cao--Shi--Zhang \cite{CSZ26}
proved H\"older continuity of the IDS for general nonconstant analytic
potentials and fixed Diophantine frequencies, including higher-dimensional
lattice models. Their quantitative exponents are determined by the
number of zeros of the sampling function minus the energy in a
complex strip, counted with multiplicity.

A major advance in the spectral theory of these operators is Avila's
global theory \cite{AvilaGlobal}, which organizes spectral energies
according to the behavior of the complexified Lyapunov exponent.
For $E\in\mathbb{R}$, let
\[
S_{E}^{v}(\theta)
=
\begin{pmatrix}
	E-v(\theta)&-1\\
	1&0
\end{pmatrix},
\qquad
S_{E,n}^{v}(\theta)
=
S_{E}^{v}(\theta+(n-1)\alpha)\cdots S_{E}^{v}(\theta).
\]

Choose $h>0$ such that $v$ extends holomorphically to the strip
$\{\theta\in\mathbb{C}/\mathbb{Z}:
|\operatorname{Im}\theta|<h\}$.
For $|\eta|<h$, the complexified Lyapunov exponent is
\[
L_{\eta}^{v}(E)
=
\lim_{n\to\infty}\frac{1}{n}
\int_{\mathbb{T}}
\log\bigl\|S_{E,n}^{v}(\theta+i\eta)\bigr\|
\,\mathrm{d}\theta,
\]
and the acceleration is defined by
\[
\omega^{v}(E)
=
\lim_{\eta\downarrow0}
\frac{L_{\eta}^{v}(E)-L_{0}^{v}(E)}{2\pi\eta}.
\]
Avila proved that the acceleration is a nonnegative integer.
An energy $E\in\Sigma_{v}$ is called supercritical if
$L_{0}^{v}(E)>0$, and subcritical if $L_{\eta}^{v}(E)=0$
for all sufficiently small $|\eta|$.
The remaining spectral energies are called critical.

The examples of Avila--Last--Shamis--Zhou \cite{ALSZ24} for the noncritical almost Mathieu operator
show that both H\"older continuity and its exponent can depend on the arithmetic properties of the
frequency. For a fixed Diophantine frequency, the IDS is locally $1/2$-H\"older throughout the
subcritical regime \cite{Av23,AJ11}. This raises the question of whether the same regularity holds
in the supercritical regime and, if not, what determines the exponent.

In his ICM~2018 survey \cite{You18}, You related this question to acceleration and conjectured
H\"older regularity of order $1/(2\omega^{v}(E_{0}))$ at a supercritical spectral energy $E_{0}$. The
conjecture was formulated for the Lyapunov exponent and is related to the corresponding IDS estimate
by the Thouless formula; see also Ge--You--Zhao \cite{GYZ}.

\Needspace{4\baselineskip}
Several results support this prediction. For trigonometric polynomial potentials of degree $k$,
Goldstein--Schlag \cite{GS08} obtained every H\"older exponent below $1/(2k)$ on intervals of
uniformly positive Lyapunov exponent, under a strong Diophantine condition. For large trigonometric
polynomial potentials, Ge--You--Zhao \cite{GYZ} obtained H\"older estimates in terms of the number
of real zeros of $v-E$, counted with multiplicity, where $E$ is a real constant. Han--Schlag
obtained the following local estimates under a strong Diophantine condition \cite{HS1} and
subsequently for all Diophantine frequencies \cite{HS2}:
\begin{equation*}
  |N_{v}(E)-N_{v}(E')|\leqslant C_{\beta} |E-E'|^{\beta}, \qquad 0<\beta<\frac{1}{2\omega^{v}(E_{0})},
\end{equation*}
for $E,E'$ near a supercritical spectral energy $E_{0}$. Ge--Jitomirskaya \cite{GJ26} proved the
endpoint $1/2$-H\"older regularity for noncritical Type~I operators. Using a transfer-operator framework for dual long-range operators,
Li--Wang--You--Zhou \cite{LWYZ26} obtained pointwise IDS estimates
with exponent $1/(2\omega^{v}(E_{0}))-\varepsilon$ for Diophantine
frequencies and sufficiently large analytic potentials, with the
coupling threshold allowed to depend on the exponent loss
$\varepsilon>0$.

The exponent $1/(2\omega^{v}(E_{0}))$ need not be sharp for a given potential. For example,
$v(\theta)=2\lambda\cos(4\pi\theta)$ with $\lambda>1$ gives an almost Mathieu operator with
frequency $2\alpha$. Its IDS is $1/2$-H\"older for Diophantine $\alpha$ \cite{AJ10}. On the
spectrum, $L_{0}^{v}(E)=\log\lambda>0$ and $\omega^{v}(E)=2$, since doubling the phase doubles the
acceleration \cite{GJYZ}.

\subsection{Main result and proof strategy}

To our knowledge, no example with an optimal local IDS H\"older exponent below $1/2$ was previously
known in the supercritical regime for analytic potentials and Diophantine frequencies. We prove that
the dependence on acceleration is sharp: for every positive integer $k$ and every Diophantine
frequency, the exponent $1/(2k)$ is attained at a supercritical spectral edge with acceleration $k$.

For an IDS $N$ and $E_{0}\in\operatorname{supp}\mathrm{d}N$, set
\begin{equation}\label{e:loc}
  \begin{aligned}
    d^{-}(N,E_{0})&=\liminf_{r\downarrow0} \frac{\log\bigl(N(E_{0}+r)-N(E_{0}-r)\bigr)}{\log(2r)},\\
    d^{+}(N,E_{0})&=\limsup_{r\downarrow0} \frac{\log\bigl(N(E_{0}+r)-N(E_{0}-r)\bigr)}{\log(2r)}.
  \end{aligned}
\end{equation}
These are the lower and upper local dimensions of $\mathrm{d}N$ at $E_{0}$. When they agree, we call
their common value the exact local H\"older exponent of $N$ at $E_{0}$; a H\"older bound at this
exponent requires a separate estimate. Throughout the paper, we assume that $\alpha$ satisfies the
Diophantine condition
\begin{equation}\label{e:dc}
  |e^{2\pi i n\alpha}-1|\geqslant\gamma |n|^{-\tau}, \qquad n\in\mathbb{Z}\setminus\{0\},
\end{equation}
for some $0<\gamma\leqslant1$ and $\tau>0$.

\begin{theorem}\label{t:main}
For every integer $k\geqslant1$ and Diophantine frequency $\alpha$, there exist a real even
trigonometric polynomial $v$ of degree $k$ and an outer edge $E_{0}$ of $\Sigma_{v}$ such that
\begin{equation*}
  L_{0}^{v}(E_{0})>0,\qquad \omega^{v}(E_{0})=k,\qquad d^{-}(N_{v},E_{0})=d^{+}(N_{v},E_{0})=\frac{1}{2k}.
\end{equation*}
\end{theorem}

For every $\beta>1/(2k)$, the IDS in Theorem~\ref{t:main} fails to be pointwise $\beta$-H\"older at
$E_{0}$. Hence the exponent cannot be improved in general, even at spectral edges. In particular,
the examples with $k\geqslant2$ are not $1/2$-H\"older there. The case $k=1$ already occurs for the
almost Mathieu operator; see Li--You--Zhou \cite{LYZ} and Aubry duality.

The construction starts from the reference polynomial
\[
V_{\mathrm{ref}}(\theta)
=(2\cos(2\pi\theta)-2)^{k},
\]
whose zero of order $2k$ gives a distribution function with local
exponent $1/(2k)$ at zero. Under Aubry duality, large coupling
corresponds to a small quasi-periodic perturbation of the
constant-coefficient finite-range operator with this Fourier symbol.
The main difficulty is to preserve the degenerate spectral edge
under this perturbation.

At zero energy, the reference dual transfer matrix is a principal
unipotent with a single Jordan block of size $2k$.
We show that varying the constant term and the $k-1$ lower Fourier
coefficients gives a family transverse to its conjugacy class.
This allows a parameter-dependent KAM iteration to analytically
conjugate the perturbed zero-energy cocycle to the same reference
matrix. The parameter correction follows the modifying-term approach
of Puig--Sim\'o \cite{PS06}; see also Hou--You \cite{HY} for reducibility of quasi-periodic linear systems.
All dual Lyapunov exponents therefore vanish at zero, and the
multiplicative Jensen formula \cite{GJYZ} gives
$L_{0}^{v}(0)=\log\lambda>0$ and $\omega^{v}(0)=k$.

To estimate the IDS, we compare the nearby-energy cocycles with the
reference family, using a perturbative argument related to those
of Moser--P\"oschel \cite{MP84} and Puig \cite{Puig06}.
The exact reduction at zero makes the comparison error
$O(\lambda^{-1}|E|)$, which is essential near the degenerate edge.
A quantitative symplectic block diagonalization then isolates a
two-dimensional center on the interior side, where a rotation-number
estimate yields matching IDS bounds of order $|E|^{1/(2k)}$.
On the exterior side, constancy of the rotation number excludes
spectrum near zero, and an operator-norm comparison establishes that
zero is an outer spectral edge. Aubry duality transfers these
estimates to the original operator and proves
Theorem~\ref{t:main}.

Section~2 introduces the dual model and collects the analytic,
symplectic and rotation-number preliminaries.
Section~3 establishes parameter transversality.
Section~4 constructs the zero-energy conjugacy and computes
the Lyapunov exponent and acceleration.
Section~5 develops the reference normal forms, proves the sharp
IDS estimates, identifies the outer spectral edge, and completes
the proof of Theorem~\ref{t:main}.

\section{Preliminaries and the dual construction}

We first fix the analytic and symplectic conventions used throughout
the paper, then introduce the dual family and formulate the quantitative
construction in Theorem~\ref{t:dual}.
The parameter selected by the KAM argument in Section~4 yields both
the dual estimates proved in Section~5 and the Lyapunov exponent and
acceleration identities needed for Theorem~\ref{t:main}.

\subsection{Analytic cocycles and Lyapunov exponents}
\label{sec:acceleration}

For $h>0$, set
\[
\mathbb{T}_{h}
=\{\theta\in\mathbb{C}/\mathbb{Z}:
|\operatorname{Im}\theta|<h\}.
\]
For a matrix-valued Fourier series, write
\[
F(\theta)
=\sum_{\ell\in\mathbb{Z}}\widehat{F}(\ell)e^{2\pi i\ell\theta},
\qquad
\|F\|_{h}
=\sum_{\ell\in\mathbb{Z}}
\|\widehat{F}(\ell)\|e^{2\pi h|\ell|},
\]
where $\|\cdot\|$ is the Euclidean operator norm.
We denote by $\mathcal{B}_{h}$ the space of series with finite
$\|\cdot\|_{h}$ norm.
For a matrix space or group $\mathcal{X}$,
$\mathcal{B}_{h}(\mathcal{X})$ denotes the functions in
$\mathcal{B}_{h}$ taking values in $\mathcal{X}$.
When $\mathcal{X}$ is real, these functions take values in
$\mathcal{X}$ on $\mathbb{T}$ and in its complexification on
$\mathbb{T}_{h}$.
We also write
\[
\langle F\rangle=\widehat{F}(0),
\qquad
\|F\|_{C^{0}}
=\sup_{\theta\in\mathbb{T}}\|F(\theta)\|.
\]

For $A\in C^{0}(\mathbb{T},\operatorname{GL}(m,\mathbb{C}))$,
the cocycle $(\alpha,A)$ is the skew product
\[
(\theta,x)\longmapsto(\theta+\alpha,A(\theta)x).
\]
Its iterates are $(\alpha,A)^{n}=(n\alpha,A^{(n)})$, where
\[
A^{(n)}(\theta)
=
\begin{cases}
	A(\theta+(n-1)\alpha)\cdots A(\theta),
	& n\geqslant1,\\[2mm]
	I_{m},
	& n=0,\\[2mm]
	A(\theta+n\alpha)^{-1}\cdots A(\theta-\alpha)^{-1},
	& n\leqslant-1.
\end{cases}
\]
An analytic cocycle is reducible on $\mathbb{T}_{h}$ if there
exist an analytic invertible matrix function $B$ on
$\mathbb{T}_{h}$ and a constant matrix
$A_{*}\in\operatorname{GL}(m,\mathbb{C})$ such that
\[
B(\theta+\alpha)^{-1}A(\theta)B(\theta)=A_{*}.
\]

The Lyapunov exponents
$L_{1}(\alpha,A)\geqslant\cdots\geqslant L_{m}(\alpha,A)$
are defined by
\[
\sum_{j=1}^{r}L_{j}(\alpha,A)
=
\lim_{n\to\infty}\frac{1}{n}
\int_{\mathbb{T}}
\log\Bigl\|
\bigwedge\nolimits^{r}A^{(n)}(\theta)
\Bigr\|\,\mathrm{d}\theta,
\qquad 1\leqslant r\leqslant m.
\]
They are invariant under continuous invertible conjugacies.
For a constant matrix, they are the logarithms of its eigenvalue
moduli, counted with algebraic multiplicity.

For $d\geqslant1$, let
\[
J_{2d}
=
\begin{pmatrix}
	0&-I_{d}\\
	I_{d}&0
\end{pmatrix}.
\]
The symplectic group is
\[
\operatorname{Sp}(2d,\mathbb{F})
=
\{M\in\operatorname{GL}(2d,\mathbb{F}):
M^{T}J_{2d}M=J_{2d}\},
\qquad \mathbb{F}\in\{\mathbb{R},\mathbb{C}\}.
\]
For cocycles taking values in
$\operatorname{Sp}(2k,\mathbb{R})$, the Lyapunov exponents occur
in opposite pairs:
\[
L_{j}(\alpha,A)=-L_{2k+1-j}(\alpha,A),
\qquad 1\leqslant j\leqslant k.
\]
In particular,
$L_{1}(\alpha,A)\geqslant\cdots\geqslant L_{k}(\alpha,A)
\geqslant0$.
For analytic real symplectic cocycles, reducing conjugacies are
also required to be symplectic and real on $\mathbb{T}$.

Throughout the paper, we write $J=J_{2k}$ and
\[
\mathfrak{g}
=\operatorname{sp}(2k,\mathbb{R})
=
\{X\in\operatorname{Mat}_{2k}(\mathbb{R}):
X^{T}J+JX=0\}.
\]
Its complexification is denoted by $\mathfrak{g}_{\mathbb{C}}$.

For symplectic matrices, $A_{1}\oplus A_{2}$ denotes the
symplectic direct sum in standard coordinates.
If $q^{(j)},p^{(j)}$ are the coordinate blocks of the $j$th
summand, this means reordering the ordinary block coordinates
by the fixed orthogonal permutation
\[
(q^{(1)},p^{(1)},q^{(2)},p^{(2)})
\longmapsto
(q^{(1)},q^{(2)},p^{(1)},p^{(2)}).
\]
We use the same coordinate order for direct sums of Lagrangian
frames.
In block calculations, we use the original order and the
corresponding block diagonal symplectic form.

We use $C$ for a positive constant that may change from line
to line.
Unless otherwise stated, constants depend only on $k$, the
Diophantine constants $\gamma,\tau$ in \eqref{e:dc}, and the
fixed analytic widths.

A special case is the Schr\"odinger cocycle
$(\alpha,S_{E}^{v})$ associated with \eqref{e:op}, where
\[
S_{E}^{v}(\theta)
=
\begin{pmatrix}
	E-v(\theta)&-1\\
	1&0
\end{pmatrix}.
\]
For $E\in\mathbb{R}$ and $|\eta|<h$, its complexified
Lyapunov exponent is
\[
L_{\eta}^{v}(E)
=
L_{1}\bigl(\alpha,S_{E}^{v}(\cdot+i\eta)\bigr),
\]
where $v$ extends holomorphically to $\mathbb{T}_{h}$
and is real on $\mathbb{T}$.
With the normalization used in the introduction, the
acceleration is
\begin{equation}\label{e:2.4}
	\omega^{v}(E)
	=
	\lim_{\eta\downarrow0}
	\frac{L_{\eta}^{v}(E)-L_{0}^{v}(E)}{2\pi\eta}.
\end{equation}
For irrational $\alpha$, the function
$\eta\mapsto L_{\eta}^{v}(E)$ is even, convex and piecewise
affine, with slopes in $2\pi\mathbb{Z}$.
In particular,
$\omega^{v}(E)\in\mathbb{Z}_{\geqslant0}$ \cite{AvilaGlobal}.

\subsection{The dual model and transfer matrices}

For $p=(a_{k-1},\ldots,a_{1},e)\in\mathbb{R}^{k}$, consider the potential
\begin{equation}\label{e:v}
	V_{p}(\theta)=2\cos(2\pi k\theta) +\sum_{j=1}^{k-1}a_{j}\cos(2\pi j\theta)+e.
\end{equation}
Let $\bar{p}$ be the unique parameter for which
\begin{equation}\label{e:v0}
	V_{\bar{p}}(\theta)=(2\cos(2\pi\theta)-2)^{k} =(-1)^{k}(2\sin(\pi\theta))^{2k}.
\end{equation}
The distribution function
$\bar{N}(E)=\operatorname{Leb}\{\theta\in\mathbb{T}:V_{\bar{p}}(\theta)\leqslant E\}$ satisfies, for
$0<r<4^{k}$,
\begin{equation*}
	\bar{N}(r)-\bar{N}(-r) =\frac{2}{\pi}\arcsin\!\biggl(\frac{r^{1/(2k)}}{2}\biggr) =\frac{1}{\pi} r^{1/(2k)}+O\!\bigl(r^{3/(2k)}\bigr).
\end{equation*}
\Needspace{5\baselineskip}
The exponent $1/(2k)$ is thus determined by the zero of order $2k$. For instance, when $k=2$,
\begin{equation*}
	V_{\bar{p}}(\theta)=2\cos(4\pi\theta)-8\cos(2\pi\theta)+6 =(2\cos(2\pi\theta)-2)^{2},
\end{equation*}
and $\bar{N}$ has local exponent $1/4$ at zero. The problem is to preserve this degeneracy at finite
coupling.

For $\lambda>0$, the Aubry dual of $\lambda^{-1}H_{\lambda V_{p},\theta}$ is the finite-range
operator
\begin{equation}\label{e:dual}
	\begin{aligned}
		(L_{p,\theta}u)_{n}={}&u_{n+k}+u_{n-k} +\sum_{j=1}^{k-1}\frac{a_{j}}{2}(u_{n+j}+u_{n-j})\\
		&+\bigl(e+2\lambda^{-1}\cos(2\pi(\theta+n\alpha))\bigr)u_{n}.
	\end{aligned}
\end{equation}
The dependence on $\lambda$ and $\alpha$ is suppressed. Denote the IDS and spectrum by
$\widetilde{N}_{p}$ and $\widetilde{\Sigma}_{p}$, respectively. By Aubry duality \cite{GYZ,HP},
\begin{equation}\label{e:au}
	\widetilde{N}_{p}(E)=N_{\lambda V_{p}}(\lambda E),\qquad \widetilde{\Sigma}_{p}=\lambda^{-1}\Sigma_{\lambda V_{p}}.
\end{equation}
At $\lambda^{-1}=0$ and $p=\bar{p}$, the Fourier transform conjugates \eqref{e:dual} to
multiplication by \eqref{e:v0}, so its IDS is $\bar{N}$. Adjusting the $k$ parameters allows us to
preserve the degenerate edge under the quasi-periodic perturbation.

\begin{theorem}\label{t:dual}
	Let $k\geqslant1$ be an integer and let $\alpha$ satisfy \eqref{e:dc}. There exist
	$\lambda_{0},c,C,E_{*}>0$, depending only on $k,\gamma,\tau$, such that for every
	$\lambda\geqslant\lambda_{0}$ there is $p_{*}\in\mathbb{R}^{k}$ with
	\begin{equation*}
		\|p_{*}-\bar{p}\|\leqslant C\lambda^{-1}.
	\end{equation*}
	For this parameter, zero is the bottom of $\widetilde{\Sigma}_{p_{*}}$ for even $k$ and the top for
	odd $k$. Moreover,
	\begin{equation*}
		cr^{1/(2k)}\leqslant\widetilde{N}_{p_{*}}(r)-\widetilde{N}_{p_{*}}(-r)\leqslant Cr^{1/(2k)},\qquad 0<r<E_{*}.
	\end{equation*}
\end{theorem}

In particular, $d^{-}(\widetilde{N}_{p_{*}},0)=d^{+}(\widetilde{N}_{p_{*}},0)=1/(2k)$, with a
pointwise H\"older bound at this exponent. The parameter $p_{*}$ may depend on $\lambda$ and
$\alpha$, which remain fixed when taking the local exponent.

For a solution of $L_{p,\theta}u=Eu$ at a real energy $E$, set
\begin{equation*}
  U_{n}=(u_{n+k-1},\ldots,u_{n},u_{n-1},\ldots,u_{n-k})^{T}, \qquad U_{n+1}=S_{E}(p,\theta+n\alpha)U_{n}.
\end{equation*}
Put $c_{k}=1$ and $c_{j}=a_{j}/2$ for $1\leqslant j<k$. The companion matrix is
\begingroup
\setlength{\arraycolsep}{4.5pt}
\begin{equation*}
  S_{E}(p,\theta)=\begin{pmatrix}
    -c_{k-1}&\cdots&-c_{1}&E-e-2\lambda^{-1}\cos(2\pi\theta) &-c_{1}&\cdots&-c_{k-1}&-1\\
    \multicolumn{7}{c}{I_{2k-1}}&0_{(2k-1)\times1}
  \end{pmatrix}.
\end{equation*}
\endgroup
The terms involving $c_{1},\ldots,c_{k-1}$ are omitted when $k=1$. Set
\begin{equation}\label{e:2.1}
  T(p)_{rs}=\begin{cases}
    c_{k+r-s},&r\leqslant s,\\
    0,&r>s,
  \end{cases} \quad 1\leqslant r,s\leqslant k, \qquad \mathcal{T}(p)=\begin{pmatrix}
    T(p)^{-1}&0\\
    0&I_{k}
  \end{pmatrix}.
\end{equation}
Thus $T(p)$ is upper triangular with diagonal entries one,
so $\mathcal{T}(p)$ is well defined.
The finite-range recurrence admits a symplectic formulation;
see \cite{HP}.
In the coordinates specified by \eqref{e:2.1}, a direct
calculation gives
\begin{equation}\label{e:2.3}
	A_{E}(p,\theta)
	=\mathcal{T}(p)^{-1}S_{E}(p,\theta)\mathcal{T}(p)
	\in\operatorname{Sp}(2k,\mathbb{R}),
	\qquad \theta\in\mathbb{T}.
\end{equation}
The symplectic form $J$ is independent of $p$. At $\lambda^{-1}=0$, we denote the
constant-coefficient operator and its transfer matrices by $\bar{L}_{p}$, $\bar{S}_{E}(p)$
and $\bar{A}_{E}(p)$. Replacing $\lambda^{-1}$ with $t\lambda^{-1}$, $0\leqslant t\leqslant1$,
joins $\bar{A}_{E}(p)$ to $A_{E}(p,\cdot)$ within $\operatorname{Sp}(2k,\mathbb{R})$. Since
this group is connected, $A_{E}(p,\cdot)$ is homotopic to the identity.

We use the multiplicative Jensen formula for real even
trigonometric polynomials
\[
v(\theta)=\sum_{j=-k}^{k}v_{j}e^{2\pi i j\theta},
\qquad
v_{-j}=v_{j}\in\mathbb{R},
\qquad v_{k}\ne0.
\]
The corresponding dual operators are
\[
(\widehat{H}_{v,\theta}u)_{n}
=
\sum_{j=-k}^{k}v_{j}u_{n+j}
+2\cos\bigl(2\pi(\theta+n\alpha)\bigr)u_{n}.
\]
Let
$\widehat{L}_{1}^{v}(E)\geqslant\cdots
\geqslant\widehat{L}_{k}^{v}(E)\geqslant0$
be the $k$ largest Lyapunov exponents of the associated
$2k$-dimensional transfer cocycle, counted with multiplicity
and normalized per scalar lattice step.
The multiplicative Jensen formula \cite{GJYZ}, together
with the duality identity \cite{HP}, gives
\begin{equation}\label{e:2.5}
	L_{\eta}^{v}(E)
	=
	\log|v_{k}|
	+\sum_{j=1}^{k}
	\max\{\widehat{L}_{j}^{v}(E),2\pi|\eta|\},
	\qquad (E,\eta)\in\mathbb{R}^{2}.
\end{equation}
Since $v$ is a trigonometric polynomial, this formula holds
for every real $\eta$.
Taking the right derivative at zero yields
\[
\omega^{v}(E)
=
\#\{1\leqslant j\leqslant k:
\widehat{L}_{j}^{v}(E)=0\}.
\]

For the family \eqref{e:dual}, we have
$\widehat{H}_{\lambda V_{p},\theta}=\lambda L_{p,\theta}$.
The equations
$\widehat{H}_{\lambda V_{p},\theta}u=\lambda Eu$
and $L_{p,\theta}u=Eu$ therefore define the same scalar
recurrence.
Consequently,
\begin{equation}\label{e:2.6}
	\widehat{L}_{j}^{\lambda V_{p}}(\lambda E)
	=
	L_{j}\bigl(\alpha,A_{E}(p,\cdot)\bigr),
	\qquad 1\leqslant j\leqslant k.
\end{equation}

\subsection{Fibered rotation numbers}
Throughout this subsection, $\alpha\in\mathbb{R}\setminus\mathbb{Q}$. We use the symplectic rotation
number \cite{Ruelle} with the Lagrangian-phase normalization of Li--Wu \cite{LW}. A real Lagrangian
frame is a matrix
\begin{equation*}
  Q=\begin{pmatrix}
    X\\
    Y
  \end{pmatrix} \in\mathbb{R}^{2k\times k}, \qquad \operatorname{rank}Q=k,\qquad Q^{T}JQ=0.
\end{equation*}
The frames $Q$ and $QR$ represent the same subspace for every $R\in\operatorname{GL}(k,\mathbb{R})$.
Since
\begin{equation*}
  (X+iY)^{*}(X+iY)=X^{T}X+Y^{T}Y=Q^{T}Q>0,
\end{equation*}
the matrices $X+iY$ and $X-iY$ are invertible. Define
\begin{equation}\label{e:2.7}
  \Phi(Q)=\det\bigl((X+iY)(X-iY)^{-1}\bigr), \qquad \phi(Q)=\frac{1}{\pi}\arg\det(X+iY)\pmod{\mathbb{Z}}.
\end{equation}
Thus $\Phi(Q)=e^{2\pi i\phi(Q)}$ and $\Phi(QR)=\Phi(Q)$. Along a path of frames, we use a continuous
real lift of $\phi$.

Let $A\in C^{0}(\mathbb{T},\operatorname{Sp}(2k,\mathbb{R}))$ be homotopic to the identity. Choose a
continuous homotopy $H_{t}(\theta)\in\operatorname{Sp}(2k,\mathbb{R})$ from $H_{0}=I$ to $H_{1}=A$.
Let $d_{A}(\theta,Q)$ be the value at $t=1$ of the continuous real lift of
$t\mapsto\Phi(H_{t}(\theta)Q)/\Phi(Q)$ that equals zero at $t=0$. Set
\begin{equation*}
  d_{A}^{(n)}(\theta,Q)=\sum_{j=0}^{n-1} d_{A}\bigl(\theta+j\alpha,A^{(j)}(\theta)Q\bigr).
\end{equation*}
These functions depend only on the subspace represented by $Q$. Any two continuous real lifts of
$\Phi(A(\theta)Q)/\Phi(Q)$ differ by an integer constant, since the
product of the torus and the real Lagrangian Grassmannian is connected. The limit
\begin{equation}\label{e:2.8}
  \rho(\alpha,A)=\lim_{n\to\infty}\frac{d_{A}^{(n)}(\theta,Q)}{n}
\end{equation}
exists uniformly in $\theta$ and the Lagrangian subspace represented by $Q$, and is independent of
both, by unique ergodicity of the irrational rotation and the rotation-number result in \cite{LW}.
Its class modulo $\mathbb{Z}$ is independent of the lift, and its local real lifts are continuous in
the $C^{0}$ topology.

If $B\in C^{0}(\mathbb{T},\operatorname{Sp}(2k,\mathbb{R}))$
is homotopic to the identity, then \cite{LW} gives
\[
\rho\bigl(
\alpha,B(\cdot+\alpha)^{-1}A(\cdot)B(\cdot)
\bigr)
=
\rho(\alpha,A)\pmod{\mathbb{Z}}.
\]
Along a homotopy $B_{t}$ from $I$ to $B$, choose continuous
real lifts of the rotation numbers of the conjugate cocycles,
starting from a fixed lift of $\rho(\alpha,A)$ at $t=0$.
Their difference from that fixed lift is continuous,
integer-valued and initially zero, so the real lifts agree.

For a constant matrix $M$, write $\rho(M)=\rho(\alpha,M)$. The path $t\mapsto R_{t\psi}$ gives
\begin{equation*}
  R_{\psi}=\begin{pmatrix}
    \cos\psi&-\sin\psi\\
    \sin\psi&\cos\psi
  \end{pmatrix}, \qquad \rho(R_{\psi})=\frac{\psi}{\pi}.
\end{equation*}
This is twice the usual oriented-vector normalization in dimension two; see \cite{LW}. With the
induced lifts, iteration and symplectic direct sums satisfy
\begin{equation}\label{e:2.10}
  \rho(n\alpha,A^{(n)})=n\rho(\alpha,A),\qquad \rho(\alpha,A_{1}\oplus A_{2})=\rho(\alpha,A_{1})+\rho(\alpha,A_{2}).
\end{equation}
The first identity follows by grouping the displacements into blocks of $n$ steps, and the second
from $\Phi(Q_{1}\oplus Q_{2})=\Phi(Q_{1})\Phi(Q_{2})$ and \eqref{e:2.8}.

\begin{lemma}\label{l:2.1}
Let $A\in C^{0}(\mathbb{T},\operatorname{Sp}(2k,\mathbb{R}))$ be homotopic to the identity. If every
$A(\theta)$ preserves the same real Lagrangian subspace, then $\rho(\alpha,A)\in\mathbb{Z}$ for
every real lift. In particular, this holds for a constant symplectic matrix preserving a real
Lagrangian subspace.
\end{lemma}

\begin{proof}
Let $Q$ be a fixed frame for the invariant subspace. For each $\theta$, write
$A(\theta)Q=QR(\theta)$ with $R(\theta)\in\operatorname{GL}(k,\mathbb{R})$. Thus
$\Phi(A(\theta)Q)=\Phi(Q)$, so $d_{A}(\theta,Q)$ is integer-valued and, by continuity, independent
of $\theta$. Every one-step displacement along the orbit equals this integer, since the subspace is
invariant. The conclusion follows from \eqref{e:2.8}.
\end{proof}

\begin{lemma}\label{l:2.2}
	\textup{(1)}
	If $K\in C^{0}(\mathbb{T},\operatorname{Sp}(2d,\mathbb{R}))$
	and $\|K-I\|_{C^{0}}\leqslant1/4$, the real lift defined by
	$t\mapsto e^{t\log K}$ satisfies
	\[
	|\rho(\alpha,K)|
	\leqslant\frac{2d}{\pi}\|K-I\|_{C^{0}}.
	\]
	
	\smallskip
	\noindent\textup{(2)}
	If $C\in C^{0}(\mathbb{T},\operatorname{SL}(2,\mathbb{R}))$
	and $\|C-R_{\psi}\|_{C^{0}}\leqslant\varepsilon\leqslant1/2$,
	the real lift whose one-step angular displacement lies in
	$(\psi-\pi/2,\psi+\pi/2)$ satisfies
	\[
	\left|\rho(\alpha,C)-\frac{\psi}{\pi}\right|
	\leqslant\frac{2\varepsilon}{\pi}.
	\]
	Both estimates are uniform in $\alpha$.
\end{lemma}

\begin{proof}
We bound the one-step phase displacement uniformly in the initial Lagrangian subspace, then
sum along the orbit and apply \eqref{e:2.8}. For a differentiable path of real Lagrangian frames
$Q(t)=\binom{X(t)}{Y(t)}$, the logarithmic derivative of the determinant in \eqref{e:2.7} gives
\begin{equation*}
  \partial_{t}\phi(Q)
  =\frac{1}{\pi}\operatorname{Im}\operatorname{tr}
  \bigl((X+iY)^{-1}(\partial_{t}X+i\partial_{t}Y)\bigr).
\end{equation*}
Since $(X+iY)^{-1}=(Q^{T}Q)^{-1}(X+iY)^{*}$ and
\begin{equation*}
  \operatorname{Im}\bigl((X+iY)^{*}(\partial_{t}X+i\partial_{t}Y)\bigr)
  =-Q^{T}J_{2d}\partial_{t}Q,
\end{equation*}
we obtain
\begin{equation*}
  \partial_{t}\phi(Q)
  =-\frac{1}{\pi}\operatorname{tr}
  \bigl((Q^{T}Q)^{-1}Q^{T}J_{2d}\partial_{t}Q\bigr).
\end{equation*}
For (1), set $Z(\theta)=\log K(\theta)$ and $Q_{t}=e^{tZ(\theta)}Q$. The logarithm is real
Hamiltonian and satisfies $\|Z\|_{C^{0}}\leqslant2\|K-I\|_{C^{0}}$. To estimate the trace
independently of the choice of frame, set $P_{t}=Q_{t}(Q_{t}^{T}Q_{t})^{-1/2}$, which has
orthonormal columns. Using $\partial_{t}Q_{t}=Z(\theta)Q_{t}$ and the cyclicity of the trace gives
\begin{equation*}
  |\partial_{t}\phi(Q_{t})|
  =\frac{1}{\pi}\bigl|\operatorname{tr}(P_{t}^{T}J_{2d}Z(\theta)P_{t})\bigr|
  \leqslant\frac{d}{\pi}\|Z(\theta)\|.
\end{equation*}
Integration in the homotopy parameter bounds each one-step phase displacement by
$2d\|K-I\|_{C^{0}}/\pi$ in absolute value. Summing over $n$ steps, dividing by $n$ and applying
\eqref{e:2.8} proves (1).

For (2), write $C(\theta)=R_{\psi}(I+G(\theta))$ with $\|G\|_{C^{0}}\leqslant\varepsilon$. The
scalar product of any unit vector $v$ with $(I+G)v$ is at least $1-\varepsilon$, while
$|\det(v,(I+G)v)|\leqslant\varepsilon$. The angle between these vectors therefore has absolute value
at most $2\varepsilon$. The lifted angular action takes the form
\begin{equation*}
  t\longmapsto t+\psi+r(\theta,t),\qquad |r(\theta,t)|\leqslant2\varepsilon.
\end{equation*}
Iteration gives $|t_{n}-t_{0}-n\psi|\leqslant2n\varepsilon$. Dividing by $\pi n$ and taking the
limit proves (2).
\end{proof}

\subsection{The IDS and rotation numbers}
We recall the IDS--rotation-number identity for real block operators \cite{LW}, with the
normalization above.

Let $\alpha\in\mathbb{R}\setminus\mathbb{Q}$, $C\in\operatorname{GL}(k,\mathbb{R})$, and let
$W:\mathbb{T}\to\mathbb{R}^{k\times k}$ be continuous and symmetric. Consider
\begin{equation*}
  (H_{\theta} u)_{n}=C^{T}u_{n-1}+W(\theta+n\alpha)u_{n}+Cu_{n+1} \quad\text{on }\ell^{2}(\mathbb{Z},\mathbb{C}^{k}).
\end{equation*}
The symplectic transfer matrix for the state $(Cu_{n},u_{n-1})^{T}$ is
\begin{equation*}
  M_{E}(\theta)=\begin{pmatrix}
    (EI_{k}-W(\theta))C^{-1}&-C^{T}\\
    C^{-1}&0
  \end{pmatrix}.
\end{equation*}
Normalize the IDS per scalar coordinate by
\begin{equation*}
  N_{b}(E)=\frac{1}{k}\int_{\mathbb{T}}\sum_{j=1}^{k} \langle\delta_{0,j},\mathbf{1}_{(-\infty,E]}(H_{\theta}) \delta_{0,j}\rangle\,\mathrm{d}\theta,
\end{equation*}
where $\delta_{0,j}$ is the $j$th coordinate vector at block zero.
Replacing $W$ with $tW$, $0\leqslant t\leqslant1$, gives a homotopy between $M_{E}$ and a constant
matrix. Since $\operatorname{Sp}(2k,\mathbb{R})$ is connected, $M_{E}$ is homotopic to the identity.

\begin{theorem}\label{t:2.3}
For every $E\in\mathbb{R}$,
\begin{equation}\label{e:2.12}
  \rho(\alpha,M_{E})=k\bigl(1-N_{b}(E)\bigr)\pmod{\mathbb{Z}}.
\end{equation}
\end{theorem}

By continuity of the IDS, every continuous real lift in \eqref{e:2.12} differs from $k(1-N_{b}(E))$
by a constant integer on $\mathbb{R}$.

For fixed $p,\lambda$, write $\widetilde{N}=\widetilde{N}_{p}$ and
$\widetilde{\Sigma}=\widetilde{\Sigma}_{p}$.

\begin{proposition}
For fixed $p,\lambda$, the cocycle \eqref{e:2.3} satisfies
\begin{equation}\label{e:2.13}
  \widetilde{N}(E)-\widetilde{N}(E') =-\bigl[\rho(\alpha,A_{E}(p,\cdot)) -\rho(\alpha,A_{E'}(p,\cdot))\bigr], \qquad E,E'\in\mathbb{R},
\end{equation}
where a continuous real lift of the rotation number is chosen as $E$ varies.
\end{proposition}

\begin{proof}
Group the scalar sites into $\Psi_{j}=(u_{jk+k-1},\ldots,u_{jk})^{T}$. For row $r$ and column $s$,
sites in adjacent blocks are separated by $k+r-s$, so the forward coupling is $C=T(p)$. The diagonal
block has entries
\begin{equation*}
  W(\theta)_{rs}=\begin{cases}
    e+2\lambda^{-1}\cos\bigl(2\pi(\theta+(k-r)\alpha)\bigr),&r=s,\\
    c_{|r-s|},&r\ne s,
  \end{cases}.
\end{equation*}
The block transfer matrix $M_{E}$ thus has $C=T(p)$ and base frequency $k\alpha$. By covariance, the
spectral measure at coordinate $r$ of block zero is the zero-site scalar measure at phase
$\theta+(k-r)\alpha$. Averaging over $\theta$ gives the same contribution $\widetilde{N}$ from each
coordinate, so the block IDS per scalar coordinate is $\widetilde{N}$. The corresponding block state
is
\begin{equation*}
  \mathcal{T}(p)^{-1}U_{jk} =\begin{pmatrix}
    T(p)\Psi_{j}\\
    \Psi_{j-1}
  \end{pmatrix}.
\end{equation*}
Thus $M_{E}=A_{E}^{(k)}$. By \eqref{e:2.10}, continuous real lifts in energy satisfy
\begin{equation*}
  \rho(k\alpha,M_{E})-\rho(k\alpha,M_{E'})=k\bigl[\rho(\alpha,A_{E})-\rho(\alpha,A_{E'})\bigr].
\end{equation*}
On the other hand, Theorem~\ref{t:2.3} gives
\begin{equation*}
  \rho(k\alpha,M_{E})-\rho(k\alpha,M_{E'})=-k\bigl[\widetilde{N}(E)-\widetilde{N}(E')\bigr].
\end{equation*}
These identities imply \eqref{e:2.13}.
\end{proof}

Finally, $\widetilde{N}$ is continuous and
$\operatorname{supp}\mathrm{d}\widetilde{N}=\widetilde{\Sigma}$; see \cite{HP}. Consequently, for
every open interval $I\subset\mathbb{R}$,
\begin{equation}\label{e:2.14}
  I\cap\widetilde{\Sigma}=\varnothing \quad\Longleftrightarrow\quad \widetilde{N}\text{ is constant on }I.
\end{equation}

\section{Parameter transversality}

The companion form of $\bar{S}_{E}(p)$ gives
\begin{equation}\label{e:3.1}
  \det(zI-\bar{A}_{E}(p))=z^{2k}+1+\sum_{j=1}^{k-1}\frac{a_{j}}{2}(z^{k+j}+z^{k-j})+(e-E)z^{k}.
\end{equation}
Set $U_{*}=\bar{A}_{0}(\bar{p})$.

\begin{proposition}\label{p:uni}
The matrix $U_{*}$ has a single Jordan block of size $2k$ at $1$, and its nilpotent logarithm $N$
satisfies
\begin{equation*}
  U_{*}=e^{N},\qquad N^{2k}=0,\qquad N^{2k-1}\ne0,\qquad N\in\mathfrak{g}.
\end{equation*}
\end{proposition}

In particular, $U_{*}$ is a principal unipotent.

\begin{proof}
Equations~\eqref{e:3.1} and \eqref{e:v0} give $\det(zI-U_{*})=(z-1)^{2k}$.

The companion matrix $S=\bar{S}_{0}(\bar{p})$ has cyclic vector $e_{1}$: for $0\leqslant r<2k$,
$S^{r}e_{1}$ has coefficient $1$ along $e_{r+1}$, and all its higher coordinates vanish. Thus
$e_{1},Se_{1},\ldots,S^{2k-1}e_{1}$ are linearly independent. The minimal polynomial of $S$, and
hence of $U_{*}$, has degree $2k$ and is therefore $(z-1)^{2k}$. This proves the first assertion.

Set $K=U_{*}-I$ and define
\begin{equation*}
  N=\sum_{j=1}^{2k-1}\frac{(-1)^{j+1}}{j} K^{j}.
\end{equation*}
The logarithm and exponential identities in $\mathbb{R}[z]/(z^{2k})$ give $e^{N}=U_{*}$, while
$N^{2k-1}=K^{2k-1}\ne0$ and $N^{2k}=0$. Each entry of $(e^{xN})^{T}Je^{xN}-J$ is a polynomial in $x$
that vanishes at every nonnegative integer, since $e^{\ell N}=U_{*}^{\ell}$ is symplectic. Thus the
matrix polynomial vanishes identically, and differentiation at $x=0$ gives $N^{T}J+JN=0$.
\end{proof}

The constant homological operator is
\begin{equation*}
  \mathcal{L}_{0}Y=U_{*}^{-1}YU_{*}-Y,\qquad Y\in\mathfrak{g}.
\end{equation*}
Equip $\mathfrak{g}$ with the Frobenius inner product $(X,Y)\mapsto\operatorname{tr}(X^{T}Y)$, and
define
\begin{equation*}
  Z=\ker\mathcal{L}_{0}^{*}=(\operatorname{Im}\mathcal{L}_{0})^{\perp}, \qquad \mathfrak{g}=\operatorname{Im}\mathcal{L}_{0}\oplus Z.
\end{equation*}
Let $\Pi_{Z}$ denote the orthogonal projection onto $Z$. The equation $\mathcal{L}_{0}Y=F$ is
solvable if and only if $\Pi_{Z}F=0$. We choose the parameters so that this condition holds. The
norm of $\Pi_{Z}$ induced by the matrix operator norm satisfies
\begin{equation*}
  \|\Pi_{Z}\|\leqslant\sqrt{2k}.
\end{equation*}
This factor is absorbed in $C(k)$. Restricting $\mathcal{L}_{0}$ to
$\mathcal{Y}=(\ker\mathcal{L}_{0})^{\perp}$ gives a unique solution for each right-hand side in
$\operatorname{Im}\mathcal{L}_{0}$, with inverse
\begin{equation*}
  R_{0}=(\mathcal{L}_{0}|_{\mathcal{Y}})^{-1}: \operatorname{Im}\mathcal{L}_{0}\longrightarrow\mathcal{Y}.
\end{equation*}
We use the same notation for the complex-linear extensions of these maps and equip them with the
norms induced by $\|\cdot\|$. The logarithmic parameter derivative
\begin{equation}\label{e:3.3}
  \overline{\mathcal{P}}\,\delta p =U_{*}^{-1}\partial_{p}\bar{A}_{0}(\bar{p})\,\delta p
\end{equation}
belongs to $\mathfrak{g}$, as follows by differentiating the symplectic identity for
$\bar{A}_{0}(p)$.

Define
\begin{equation*}
  \overline{\mathcal{M}}:\mathbb{R}^{k}\longrightarrow Z,\qquad \overline{\mathcal{M}}\,\delta p=\Pi_{Z}\overline{\mathcal{P}}\,\delta p.
\end{equation*}
\begin{proposition}\label{p:tr}
We have $\dim Z=k$, and $\overline{\mathcal{M}}$ is an isomorphism with
$\|\overline{\mathcal{M}}^{-1}\|\leqslant C(k)$.
\end{proposition}

\begin{proof}
Since $U_{*}$ and $N$ are polynomials in each other, a matrix commutes with $U_{*}$ if and only if
it commutes with $N$. Choose $v$ with $N^{2k-1}v\ne0$. Then $v,Nv,\ldots,N^{2k-1}v$ form a basis. If
$Y$ commutes with $N$, write $Yv=f(N)v$ with $\deg f<2k$; then $YN^{r}v=f(N)N^{r}v$ for every basis
vector, and hence $Y=f(N)$. Conversely, every polynomial in $N$ commutes with it. Since
\begin{equation*}
  f(N)^{T}J+Jf(N)=J\bigl(f(-N)+f(N)\bigr)
\end{equation*}
and $I,N,\ldots,N^{2k-1}$ are linearly independent, such a polynomial belongs to $\mathfrak{g}$ if
and only if its even coefficients vanish. Hence
\begin{equation*}
  \ker\mathcal{L}_{0}=\operatorname{span}\{N,N^{3},\ldots,N^{2k-1}\}, \qquad \dim Z=\dim\ker\mathcal{L}_{0}=k.
\end{equation*}

Suppose now that $\overline{\mathcal{M}}\delta p=0$. Since
$\ker\Pi_{Z}=\operatorname{Im}\mathcal{L}_{0}$, there exists $Y\in\mathfrak{g}$ such that
$\overline{\mathcal{P}}\delta p=\mathcal{L}_{0}Y$. Equation~\eqref{e:3.3} then gives
\begin{equation*}
  \partial_{p}\bar{A}_{0}(\bar{p})\,\delta p
  =YU_{*}-U_{*}Y
  =\left.\frac{\mathrm{d}}{\mathrm{d}t}
    \bigl(e^{tY}U_{*}e^{-tY}\bigr)\right|_{t=0}.
\end{equation*}
Thus the parameter variation is tangent to the conjugacy class of $U_{*}$. Conjugacy invariance of
the characteristic polynomial implies that its derivative in this direction vanishes. Writing
$\delta p=(\delta a_{k-1},\ldots,\delta a_{1},\delta e)$, we obtain from \eqref{e:3.1} the
polynomial identity
\begin{equation*}
  \begin{aligned}
    0
    &=\left.\frac{\mathrm{d}}{\mathrm{d}t}
      \det\bigl(zI-\bar{A}_{0}(\bar{p}+t\delta p)\bigr)\right|_{t=0}\\
    &=\sum_{j=1}^{k-1}\frac{\delta a_{j}}{2}
      (z^{k+j}+z^{k-j})+\delta e\,z^{k}.
  \end{aligned}
\end{equation*}
The coefficients of $z^{k+j}$, $1\leqslant j<k$, and $z^{k}$ give $\delta a_{j}=0$ and $\delta e=0$,
respectively. Thus $\overline{\mathcal{M}}$ is injective and, since $\dim Z=k$, is an isomorphism.
The inverse bound follows since the map and the chosen norms depend only on $k$.
\end{proof}

\section{Reducibility at zero energy}

Fix $0<h_{*}<h_{0}$. We choose the parameter so that the cocycle at zero energy is reducible to
$U_{*}$. The KAM scheme is related to the parameter-dependent reducibility method of Puig--Sim\'o
\cite{PS06} for quasi-periodic differential equations. Here Proposition~\ref{p:tr} allows us to
eliminate the constant obstruction by varying the coefficients of $V_{p}$.

\begin{proposition}\label{p:4.1}
There exist $\lambda_{0},C>0$, depending only on $k,\gamma,\tau,h_{0},h_{*}$, such that, for every
$\lambda\geqslant\lambda_{0}$, there are $p_{*}\in\mathbb{R}^{k}$ and
$B_{*}\in\mathcal{B}_{h_{*}}(\operatorname{Sp}(2k,\mathbb{R}))$ satisfying
\begin{equation}\label{e:4.1}
  B_{*}(\theta+\alpha)^{-1}A_{0}(p_{*},\theta)B_{*}(\theta)=U_{*}.
\end{equation}
Moreover,
\begin{equation}\label{e:4.2}
  \|p_{*}-\bar{p}\|\leqslant C\lambda^{-1},\qquad \|B_{*}-I\|_{h_{*}}\leqslant C\lambda^{-1}.
\end{equation}
\end{proposition}

On $\mathfrak{g}_{\mathbb{C}}$-valued analytic functions, define
\begin{equation*}
  (\mathcal{L}X)(\theta)=U_{*}^{-1}X(\theta+\alpha)U_{*}-X(\theta).
\end{equation*}
Its restriction to constants is $\mathcal{L}_{0}$.

\begin{lemma}\label{l:4.2}
Let $0<h'<h\leqslant h_{0}$ and $F\in\mathcal{B}_{h}(\mathfrak{g})$ satisfy
$\Pi_{Z}\langle F\rangle=0$. Then $\mathcal{L}X=F$ has a unique solution
$X\in\mathcal{B}_{h'}(\mathfrak{g})$ with $\langle X\rangle\in\mathcal{Y}$ and
\begin{equation}\label{e:4.3}
  \|X\|_{h'}\leqslant \bigl(\|R_{0}\|+C_{1}\gamma^{-(4k-1)} (h-h')^{-(4k-1)\tau}\bigr)\|F\|_{h},
\end{equation}
where $C_{1}$ depends only on $k,\tau$. If $\langle F\rangle=0$, then $\langle X\rangle=0$ and the
term $\|R_{0}\|$ can be omitted.
\end{lemma}

\begin{proof}
By Proposition~\ref{p:uni}, $U_{*}=e^{N}$ with $N^{2k}=0$. Set $\operatorname{ad}_{N}(W)=NW-WN$. The
identity
\begin{equation*}
  (\operatorname{ad}_{N})^{r}W=\sum_{j=0}^{r}(-1)^{j}\binom{r}{j} N^{r-j}WN^{j}
\end{equation*}
together with $N^{2k}=0$ gives $(\operatorname{ad}_{N})^{4k-1}=0$. Since
$\mathcal{L}_{0}=e^{-\operatorname{ad}_{N}}-I$ is a polynomial in $\operatorname{ad}_{N}$ with zero
constant term,
\begin{equation*}
  \mathcal{L}_{0}^{4k-1}=0.
\end{equation*}
For $\ell\ne0$, put $\zeta_{\ell}=e^{2\pi i\ell\alpha}$. The operator on the $\ell$th Fourier mode
is $\mathcal{L}_{\ell}=(\zeta_{\ell}-1)I+\zeta_{\ell}\mathcal{L}_{0}$, whose inverse is the finite
sum
\begin{equation*}
  \mathcal{L}_{\ell}^{-1} =\sum_{j=0}^{4k-2} \frac{(-\zeta_{\ell}\mathcal{L}_{0})^{j}}{(\zeta_{\ell}-1)^{j+1}}.
\end{equation*}
The Diophantine condition \eqref{e:dc} and $|\zeta_{\ell}-1|\leqslant2$ yield
\begin{equation*}
  \|\mathcal{L}_{\ell}^{-1}\| \leqslant C\gamma^{-(4k-1)}|\ell|^{(4k-1)\tau}.
\end{equation*}
Since the zero-mode condition is equivalent to
$\langle F\rangle\in\operatorname{Im}\mathcal{L}_{0}$, the normalized solution is given by
\begin{equation*}
  \widehat{X}(0)=R_{0}\langle F\rangle,\qquad \widehat{X}(\ell)=\mathcal{L}_{\ell}^{-1}\widehat{F}(\ell) \quad(\ell\ne0).
\end{equation*}
The nonzero modes satisfy
\begin{equation*}
  \begin{aligned}
    \sum_{\ell\ne0}\|\widehat{X}(\ell)\|e^{2\pi h'|\ell|} &\leqslant C\gamma^{-(4k-1)} \sup_{x>0}\bigl(x^{(4k-1)\tau}e^{-2\pi(h-h')x}\bigr) \|F-\langle F\rangle\|_{h}\\
    &\leqslant C_{1}\gamma^{-(4k-1)}(h-h')^{-(4k-1)\tau} \|F-\langle F\rangle\|_{h}.
  \end{aligned}
\end{equation*}
Together with the zero-mode bound, this proves \eqref{e:4.3}. All operators preserve
$\mathfrak{g}_{\mathbb{C}}$, and $\mathcal{L}_{-\ell}^{-1}=\overline{\mathcal{L}_{\ell}^{-1}}$, so
the solution is real on the real torus. Uniqueness follows from invertibility on the nonzero modes
and injectivity of $\mathcal{L}_{0}$ on $\mathcal{Y}$.
\end{proof}

Define the logarithmic parameter derivative by
\begin{equation*}
  \mathcal{P}(p,\theta)=A_{0}(p,\theta)^{-1}\partial_{p}A_{0}(p,\theta),
\end{equation*}
viewed as a linear map from parameter increments to $\mathfrak{g}_{\mathbb{C}}$. For such functions,
$\|\cdot\|_{h}$ denotes the weighted Fourier norm induced by the operator norm from the Euclidean
space $\mathbb{C}^{k}$ to $\mathfrak{g}_{\mathbb{C}}$ with the matrix norm $\|\cdot\|$. The map
$\mathcal{P}$ evaluated at $p=\bar{p}$ and $\lambda^{-1}=0$ is $\overline{\mathcal{P}}$ in
\eqref{e:3.3}.

\begin{lemma}\label{l:4.3}
There exist $\lambda_{1}\geqslant1$, $r>0$, $0<b<1/4$, $C_{2}\geqslant1$ and $\eta_{*}>0$, depending
only on $k,h_{0}$, such that the following holds. Let $\lambda\geqslant\lambda_{1}$,
$0<h'<h\leqslant h_{0}$, and
\begin{equation*}
  K\geqslant\max\bigl\{1,\|R_{0}\|+ C_{1}\gamma^{-(4k-1)}(h-h')^{-(4k-1)\tau}\bigr\}.
\end{equation*}
Suppose $p_{n}\in\mathbb{R}^{k}$, $B_{n}\in\mathcal{B}_{h}(\operatorname{Sp}(2k,\mathbb{R}))$ and
$F_{n}\in\mathcal{B}_{h}(\mathfrak{g})$ satisfy
\begin{equation*}
  B_{n}(\theta+\alpha)^{-1}A_{0}(p_{n},\theta)B_{n}(\theta) =U_{*}e^{F_{n}(\theta)},\qquad \|p_{n}-\bar{p}\|\leqslant r,\qquad \|B_{n}-I\|_{h}\leqslant b.
\end{equation*}
If $f_{n}=\|F_{n}\|_{h}$ satisfies $Kf_{n}\leqslant\eta_{*}$, there exist
$\Delta p_{n}\in\mathbb{R}^{k}$ and $X_{n},F_{n+1}\in\mathcal{B}_{h'}(\mathfrak{g})$ such that, with
$p_{n+1}=p_{n}+\Delta p_{n}$ and $B_{n+1}=B_{n}e^{X_{n}}$,
\begin{equation*}
  B_{n+1}(\theta+\alpha)^{-1}A_{0}(p_{n+1},\theta)B_{n+1}(\theta) =U_{*}e^{F_{n+1}(\theta)}.
\end{equation*}
They satisfy
\begin{equation}\label{e:4.4}
  \|\Delta p_{n}\|\leqslant C_{2}f_{n},\qquad \|X_{n}\|_{h'}\leqslant C_{2}Kf_{n},
\end{equation}
and
\begin{equation}\label{e:4.5}
  \|F_{n+1}\|_{h'}\leqslant C_{2}K^{2}f_{n}^{2}.
\end{equation}
\end{lemma}

\begin{proof}
The matrix $T(p)$ is unit upper triangular, so $T(p)^{-1}$ is polynomial in $p$. Consequently,
$A_{0}$, its inverse and its first two parameter derivatives are uniformly bounded in
$\mathcal{B}_{h_{0}}$ for $\|p-\bar{p}\|\leqslant2$ and $\lambda\geqslant1$, with bounds depending
only on $k,h_{0}$. Analyticity in $p$ and $\lambda^{-1}$ gives
\begin{equation*}
  \|\mathcal{P}(p)-\overline{\mathcal{P}}\|_{h_{0}} \leqslant C\bigl(\lambda^{-1}+\|p-\bar{p}\|\bigr).
\end{equation*}

Define
\begin{equation*}
  \mathcal{A}_{n}(p,\theta) =B_{n}(\theta+\alpha)^{-1}A_{0}(p,\theta)B_{n}(\theta).
\end{equation*}
Keeping $B_{n}$ fixed and differentiating with respect to $p$ gives
\begin{equation*}
  \mathcal{P}_{n} :=\mathcal{A}_{n}(p_{n})^{-1}\partial_{p}\mathcal{A}_{n}(p_{n}) =B_{n}^{-1}\mathcal{P}(p_{n})B_{n}.
\end{equation*}
If $\|B_{n}-I\|_{h}\leqslant1/2$, the Neumann series bounds $B_{n}^{-1}$ uniformly, and hence
$\|\mathcal{P}_{n}\|_{h}\leqslant C$. Set $\mathcal{M}_{n}=\Pi_{Z}\langle\mathcal{P}_{n}\rangle$.
Then
\begin{equation*}
  \|\mathcal{M}_{n}-\overline{\mathcal{M}}\| \leqslant C\bigl(\lambda^{-1}+\|p_{n}-\bar{p}\| +\|B_{n}-I\|_{h}\bigr).
\end{equation*}
Here $C$ is uniform on the fixed neighborhoods $\|p-\bar{p}\|\leqslant1$ and
$\|B-I\|_{h}\leqslant1/2$. By Proposition~\ref{p:tr}, we may choose $0<r\leqslant1$, $0<b<1/4$ and
$\lambda_{1}\geqslant1$ so that
\begin{equation*}
  C\|\overline{\mathcal{M}}^{-1}\| (\lambda_{1}^{-1}+r+b)\leqslant\frac{1}{2}.
\end{equation*}
The Neumann series then gives
\begin{equation}\label{e:4.6}
  \|\mathcal{M}_{n}^{-1}\|\leqslant2\|\overline{\mathcal{M}}^{-1}\|.
\end{equation}

For the updates $p_{n+1}=p_{n}+\Delta p_{n}$ and $B_{n+1}=B_{n}e^{X_{n}}$, expansion of the new
error to first order jointly in $F_{n}$, $X_{n}$ and $\Delta p_{n}$ gives
$F_{n}+\mathcal{P}_{n}\Delta p_{n}-\mathcal{L}X_{n}$. The solvability condition
$\Pi_{Z}\langle F_{n}+\mathcal{P}_{n}\Delta p_{n}\rangle=0$ determines
\begin{equation*}
  \Delta p_{n}=-\mathcal{M}_{n}^{-1}\Pi_{Z}\langle F_{n}\rangle.
\end{equation*}
This choice gives $\|\Delta p_{n}\|\leqslant Cf_{n}$ and
$\|F_{n}+\mathcal{P}_{n}\Delta p_{n}\|_{h}\leqslant Cf_{n}$. Lemma~\ref{l:4.2} gives the unique
normalized solution of
\begin{equation}\label{e:4.7}
  \mathcal{L}X_{n}=F_{n}+\mathcal{P}_{n}\Delta p_{n},\qquad \langle X_{n}\rangle\in\mathcal{Y},
\end{equation}
with $\|X_{n}\|_{h'}\leqslant CKf_{n}$.

There is $0<d_{*}\leqslant1$, uniform on these neighborhoods, such that
$\mathcal{A}_{n}(p_{n})^{-1}\mathcal{A}_{n}(p_{n}+d)$ has a local logarithm for
$\|d\|\leqslant d_{*}$. The factors $B_{n}(\theta+\alpha)$ cancel in this product. Taylor's formula
gives
\begin{equation*}
  \mathcal{A}_{n}(p_{n})^{-1}\mathcal{A}_{n}(p_{n}+d) =I+\mathcal{P}_{n}d+O(\|d\|^{2})
\end{equation*}
with a uniform remainder in $\mathcal{B}_{h}$, since the parameter segment remains in
$\|p-\bar{p}\|\leqslant2$. Thus
\begin{equation*}
  \begin{aligned}
    \mathcal{A}_{n}(p_{n}+d) &=\mathcal{A}_{n}(p_{n})e^{\mathcal{P}_{n}d+\mathcal{R}_{n}(d)},\\
    \|\mathcal{R}_{n}(d)\|_{h}&\leqslant C\|d\|^{2},\qquad \|d\|\leqslant d_{*}.
  \end{aligned}
\end{equation*}
For real $d$, both the logarithm and $\mathcal{R}_{n}(d)$ belong to $\mathcal{B}_{h}(\mathfrak{g})$.
This follows from $\log(M^{-1})=-\log M$ and $M^{-1}=J^{-1}M^{T}J$ for symplectic $M$ near $I$.

Put $V_{n}=U_{*}^{-1}X_{n}(\cdot+\alpha)U_{*}$ and write
\begin{equation*}
  B_{n+1}(\theta+\alpha)^{-1}A_{0}(p_{n+1},\theta)B_{n+1}(\theta)=U_{*}\mathcal{E}_{n}(\theta).
\end{equation*}
The conjugacy identity gives
\begin{equation*}
  \mathcal{E}_{n} =e^{-V_{n}}e^{F_{n}} e^{\mathcal{P}_{n}\Delta p_{n}+\mathcal{R}_{n}(\Delta p_{n})}e^{X_{n}}.
\end{equation*}
The linear terms cancel by \eqref{e:4.7}:
\begin{equation*}
  -V_{n}+F_{n}+\mathcal{P}_{n}\Delta p_{n}+X_{n}=0.
\end{equation*}
Each of the four exponents has $\mathcal{B}_{h'}$ norm at most $CKf_{n}$ provided $Kf_{n}$ is
sufficiently small. Moreover, $\|\mathcal{R}_{n}(\Delta p_{n})\|_{h'}\leqslant Cf_{n}^{2}$.
Expanding the exponentials as $I+Z+O(\|Z\|^{2})$ in the Banach algebra and using the cancellation
above yields
\begin{equation*}
  \|\mathcal{E}_{n}-I\|_{h'}\leqslant CK^{2}f_{n}^{2}.
\end{equation*}

All constants above, including that in \eqref{e:4.6}, depend only on $k,h_{0}$. Using the
exponential estimates on the unit ball, choose $C_{2}$ to dominate the constants in the correction
and exponent bounds, and twice the constant in the bound for $\mathcal{E}_{n}-I$. Choose
\begin{equation*}
  0<\eta_{*}\leqslant \min\biggl\{\frac{d_{*}}{C_{2}},\frac{1}{C_{2}}, \frac{1}{2\sqrt{C_{2}}}\biggr\}.
\end{equation*}
When $Kf_{n}\leqslant\eta_{*}$, the parameter segment $p_{n}+t\Delta p_{n}$,
$0\leqslant t\leqslant1$, remains in $\|p-\bar{p}\|\leqslant2$, every exponent has norm at most one,
and $\|\mathcal{E}_{n}-I\|_{h'}\leqslant1/4$. Hence $F_{n+1}=\log\mathcal{E}_{n}$ is well defined
and satisfies \eqref{e:4.5}, since $\|\log(I+Z)\|\leqslant2\|Z\|$ for $\|Z\|\leqslant1/2$. The local
logarithm argument shows that it belongs to $\mathcal{B}_{h'}(\mathfrak{g})$, and
$B_{n+1}=B_{n}e^{X_{n}}$ is symplectic.
\end{proof}

\begin{proof}[Proof of Proposition~\ref{p:4.1}]
Set
\begin{equation*}
  h_{n}=h_{*}+(h_{0}-h_{*})2^{-n}.
\end{equation*}
Lemma~\ref{l:4.2} allows us to choose $K\geqslant1$, depending only on $k,\gamma,\tau,h_{0},h_{*}$,
so that
\begin{equation*}
  \|R_{0}\|+C_{1}\gamma^{-(4k-1)}(h_{n}-h_{n+1})^{-(4k-1)\tau} \leqslant K_{n}:=K2^{(4k-1)\tau n}.
\end{equation*}
Choose $C_{3}\geqslant1$, depending only on $k,h_{0}$, such that, for all sufficiently large
$\lambda$,
\begin{equation*}
  F_{0}=\log\bigl(U_{*}^{-1}A_{0}(\bar{p})\bigr),\qquad \|F_{0}\|_{h_{0}}\leqslant C_{3}\lambda^{-1}.
\end{equation*}
The estimate follows from $\|U_{*}^{-1}A_{0}(\bar{p})-I\|_{h_{0}}=O(\lambda^{-1})$, and the local
symplectic logarithm gives $F_{0}\in\mathcal{B}_{h_{0}}(\mathfrak{g})$.

Fix the constants of Lemma~\ref{l:4.3} and choose $0<\bar{\varepsilon}<1$ so that
\begin{equation}\label{e:4.8}
  \begin{gathered}
    2^{2(4k-1)\tau}\bar{\varepsilon}^{1/4}\leqslant\frac{1}{2},\qquad C_{2}K^{2}\bar{\varepsilon}^{1/2}\leqslant1,\\
    K\bar{\varepsilon}\leqslant \min\biggl\{\eta_{*},\frac{r}{2C_{2}}, \frac{\log(1+b)}{2C_{2}}\biggr\}.
  \end{gathered}
\end{equation}
Take $\lambda_{0}\geqslant\max\{\lambda_{1},C_{3}/\bar{\varepsilon}\}$ sufficiently large for the
initial logarithm estimate. For each $\lambda\geqslant\lambda_{0}$, set
\begin{equation*}
  \varepsilon_{0}=C_{3}\lambda^{-1},\qquad \varepsilon_{n}=\varepsilon_{0}^{(3/2)^{n}}.
\end{equation*}
The first two conditions in \eqref{e:4.8} imply
\begin{equation*}
  C_{2}K_{n}^{2}\varepsilon_{n}^{1/2}\leqslant1,\qquad \frac{\varepsilon_{n+1}}{\varepsilon_{n}}\leqslant\frac{1}{2},\qquad \frac{K_{n+1}\varepsilon_{n+1}}{K_{n}\varepsilon_{n}}\leqslant\frac{1}{2}.
\end{equation*}
For the first bound, the ratio of consecutive left-hand sides is
$2^{2(4k-1)\tau}\varepsilon_{n}^{1/4}\leqslant1/2$. Consequently,
\begin{equation}\label{e:4.9}
  \sum_{j=n}^{\infty}\varepsilon_{j}\leqslant2\varepsilon_{n},\qquad \sum_{j=n}^{\infty}K_{j}\varepsilon_{j}\leqslant2K_{n}\varepsilon_{n},\qquad K_{n}\varepsilon_{n}\leqslant K\varepsilon_{0}\leqslant\eta_{*}.
\end{equation}

Starting from $p_{0}=\bar{p}$, $B_{0}=I$ and the above $F_{0}$, we apply Lemma~\ref{l:4.3} with
$h=h_{n}$, $h'=h_{n+1}$ and $K=K_{n}$ and prove by induction that
\begin{equation}\label{e:4.10}
  B_{n}(\theta+\alpha)^{-1}A_{0}(p_{n},\theta)B_{n}(\theta) =U_{*}e^{F_{n}(\theta)},
\end{equation}
and
\begin{equation*}
  \|F_{n}\|_{h_{n}}\leqslant\varepsilon_{n},\qquad \|p_{n}-\bar{p}\|\leqslant r,\qquad \|B_{n}-I\|_{h_{n}}\leqslant b.
\end{equation*}
The smallness condition holds since $K_{n}\|F_{n}\|_{h_{n}}\leqslant\eta_{*}$. The quadratic
estimate gives
\begin{equation*}
  \|F_{n+1}\|_{h_{n+1}} \leqslant C_{2}K_{n}^{2}\varepsilon_{n}^{2} \leqslant\varepsilon_{n}^{3/2}=\varepsilon_{n+1}.
\end{equation*}
The parameter bound follows from \eqref{e:4.4} and \eqref{e:4.9}:
\begin{equation*}
  \|p_{n+1}-\bar{p}\| \leqslant C_{2}\sum_{j=0}^{n}\varepsilon_{j} \leqslant2C_{2}\varepsilon_{0}\leqslant r.
\end{equation*}
Since $B_{n+1}=e^{X_{0}}\cdots e^{X_{n}}$ and each earlier $X_{j}$ is defined on
$\mathbb{T}_{h_{n+1}}$, we have
\begin{equation*}
  \begin{aligned}
    \|B_{n+1}-I\|_{h_{n+1}} &\leqslant\exp\biggl(\sum_{j=0}^{n}\|X_{j}\|_{h_{n+1}}\biggr)-1\\
    &\leqslant\exp(2C_{2}K\varepsilon_{0})-1\leqslant b.
  \end{aligned}
\end{equation*}
Thus the parameter and conjugacy remain in the prescribed neighborhoods, completing the induction
for every $n$.

The tail estimates imply $p_{n}\to p_{*}$ in $\mathbb{R}^{k}$ and $B_{n}\to B_{*}$ in
$\mathcal{B}_{h_{*}}$, since
\begin{equation*}
  \|p_{n+1}-p_{n}\|\leqslant C_{2}\varepsilon_{n},\qquad \|B_{n+1}-B_{n}\|_{h_{*}}\leqslant CK_{n}\varepsilon_{n}.
\end{equation*}
The bound $\|B_{n}-I\|_{h_{*}}\leqslant b<1/4$ also ensures convergence of the inverses. The limit
is real and symplectic. Since $\|F_{n}\|_{h_{*}}\to0$, passing to the limit in \eqref{e:4.10} yields
\eqref{e:4.1}. The limiting estimates are
\begin{equation*}
  \begin{aligned}
    \|p_{*}-\bar{p}\|&\leqslant2C_{2}\varepsilon_{0},\\
    \|B_{*}-I\|_{h_{*}} &\leqslant e^{2C_{2}K\varepsilon_{0}}-1 \leqslant2(1+b)C_{2}K\varepsilon_{0}.
  \end{aligned}
\end{equation*}
Substituting $\varepsilon_{0}=C_{3}\lambda^{-1}$ proves \eqref{e:4.2}, with the stated dependence of
the constants.
\end{proof}

\begin{corollary}\label{c:4.4}
For $\lambda\geqslant\lambda_{0}$, let $p_{*}$ be given by Proposition~\ref{p:4.1} and
$v=\lambda V_{p_{*}}$. Then
\begin{equation}\label{e:4.11}
  L_{\eta}^{v}(0)=\log\lambda+2\pi k|\eta|, \qquad \eta\in\mathbb{R}.
\end{equation}
In particular,
\begin{equation}\label{e:4.12}
  L_{0}^{v}(0)=\log\lambda,\qquad \omega^{v}(0)=k.
\end{equation}
\end{corollary}

\begin{proof}
Equations~\eqref{e:2.6} and \eqref{e:4.1}, together with invariance of Lyapunov exponents under
continuous conjugacies, give
\begin{equation*}
  \widehat{L}_{j}^{v}(0) =L_{j}\bigl(\alpha,A_{0}(p_{*},\cdot)\bigr) =L_{j}(\alpha,U_{*})=0,\qquad 1\leqslant j\leqslant k.
\end{equation*}
The last equality follows from $U_{*}=e^{N}$ with $N$ nilpotent, since all eigenvalues of $U_{*}$
equal one. By \eqref{e:v}, the highest Fourier coefficient of $v$ is $v_{k}=\lambda$. The
multiplicative Jensen formula \eqref{e:2.5} gives \eqref{e:4.11} for every $\eta\in\mathbb{R}$.
Evaluating at $\eta=0$ and taking the right derivative as in \eqref{e:2.4} gives \eqref{e:4.12}.
\end{proof}

\section{Proof of Theorem~\ref{t:dual}}

Throughout this section, norms of matrix functions are $C^{0}$ norms on $\mathbb{T}$. In pointwise
identities, omitted arguments are understood to be $\theta$.

\subsection{Normal form of the reference family}

We first put the reference transfer matrices associated with
\eqref{e:v0} in symplectic normal form.
The following estimates control the change of coordinates
near the degenerate edge.

\begin{lemma}[Reference symplectic coordinates]\label{l:short-reference}
	There exist $E_{*},c,C>0$, depending only on $k$, such that for $0<|E|<E_{*}$ and $s=|E|^{1/(2k)}$
	there is a constant real symplectic matrix $W_{E}$ satisfying
	\begin{equation}\label{e:short-condition}
		\|W_{E}\|+\|W_{E}^{-1}\|\leqslant Cs^{-(2k-1)/2}.
	\end{equation}
	In these coordinates,
	\begin{equation*}
		\bar{M}_{E}:=W_{E}^{-1}\bar{A}_{E}(\bar{p})W_{E}
		=\begin{cases}
			R_{\psi_{E}}\oplus\operatorname{diag}(D_{E},D_{E}^{-T}),&(-1)^{k}E>0,\\
			\operatorname{diag}(D_{E},D_{E}^{-T}),&(-1)^{k}E<0,
		\end{cases}
	\end{equation*}
	where $D_{E}$ has dimension $k-1$ in the first case and $k$ in the second, and
	$|\psi_{E}|=2\arcsin(s/2)$. Moreover,
	\begin{equation}\label{e:short-model-gap}
		\|D_{E}\|\leqslant e^{-cs},\qquad \|\bar{M}_{E}-I\|\leqslant Cs.
	\end{equation}
	Empty hyperbolic blocks are omitted.
\end{lemma}

\begin{proof}
	The companion form and \eqref{e:v0} give
	\begin{equation*}
		\det(zI-\bar{A}_{E}(\bar{p}))
		=z^{k}\bigl((z+z^{-1}-2)^{k}-E\bigr)
		=(z-1)^{2k}-Ez^{k}.
	\end{equation*}
	Fix the sign $\sigma=E/|E|$ and write $E=\sigma s^{2k}$. Substituting $z=1+s\mu$ in the
	characteristic equation gives
	\begin{equation*}
		\mu^{2k}=\sigma(1+s\mu)^{k}.
	\end{equation*}
	At $s=0$, this equation has $2k$ distinct roots $\mu_{1},\ldots,\mu_{2k}$, with
	$\mu_{j}^{2k}=\sigma$. The derivative of $\mu^{2k}-\sigma(1+s\mu)^{k}$ with respect to $\mu$ at
	$(s,\mu)=(0,\mu_{j})$ is $2k\mu_{j}^{2k-1}\ne0$. The analytic implicit function theorem therefore
	gives an analytic root $\mu_{j}(s)$ near each $\mu_{j}$, with $\mu_{j}(0)=\mu_{j}$. Hence
	\begin{equation*}
		z_{j}(s)=1+s\mu_{j}(s)=1+s\mu_{j}+O(s^{2}),\qquad 1\leqslant j\leqslant2k.
	\end{equation*}
	For sufficiently small $s>0$, these roots are distinct and give all the eigenvalues of
	$\bar{A}_{E}(\bar{p})$.
	
	The expansion implies
	\begin{equation*}
		\log|z_{j}(s)|=s\operatorname{Re}\mu_{j}+O(s^{2}).
	\end{equation*}
	Thus, whenever $\operatorname{Re}\mu_{j}\ne0$, after decreasing $E_{*}$ we have
	\begin{equation*}
		|\log|z_{j}(s)||\geqslant cs,
	\end{equation*}
	where $c>0$ depends only on $k$. Indeed, for the two choices of $\sigma$, there are only finitely
	many nonzero values of $|\operatorname{Re}\mu_{j}|$. Moreover, $|z_{j}(s)|<1$ if
	$\operatorname{Re}\mu_{j}<0$, and $|z_{j}(s)|>1$ if $\operatorname{Re}\mu_{j}>0$.
	
	All roots of $\mu^{2k}=\sigma$ have modulus one, so the only possible purely imaginary roots are
	$\mu=\pm i$. These occur precisely when $\sigma=(-1)^{k}$. In this case, the characteristic
	equation has the two exact roots $e^{\pm2i\arcsin(s/2)}$, since
	\begin{equation*}
		(e^{it}+e^{-it}-2)^{k}=(-1)^{k}(2\sin(t/2))^{2k}.
	\end{equation*}
	Their expansions are $1\pm is+O(s^{2})$, so they are the roots corresponding to $\mu=\pm i$.
	Every other root satisfies the preceding hyperbolic estimate. If $\sigma\ne(-1)^{k}$, no
	$\mu_{j}$ is purely imaginary, and the estimate holds for every eigenvalue. The roots with
	nonzero real part occur in pairs $\mu,-\mu$. Consequently, there are $k-1$ stable and $k-1$
	unstable eigenvalues when $(-1)^{k}E>0$, and $k$ of each when $(-1)^{k}E<0$.
	
	For each root $z$, choose the eigenvector
	\begin{equation*}
		\widehat{u}(z)=\mathcal{T}(\bar{p})^{-1}(z^{2k-1},z^{2k-2},\ldots,z,1)^{T}.
	\end{equation*}
	We have $\|\widehat{u}(z)\|\leqslant C$.
	A direct computation from the definition of $\mathcal{T}(\bar{p})$ gives
	\begin{equation*}
		\begin{aligned}
			\widehat{u}(z^{-1})^{T}J\widehat{u}(z)
			&=\sum_{\ell=1}^{k}\ell c_{\ell}(z^{\ell}-z^{-\ell})\\
			&=z\frac{\mathrm{d}}{\mathrm{d}z}(z+z^{-1}-2)^{k}\\
			&=kz^{-k}(z+1)(z-1)^{2k-1},
		\end{aligned}
	\end{equation*}
	where the coefficients $c_{\ell}$ are evaluated at $p=\bar{p}$. Since
	$z_{j}(s)-1=s\mu_{j}+O(s^{2})$, this implies
	\begin{equation*}
		cs^{2k-1}\leqslant|\widehat{u}(z^{-1})^{T}J\widehat{u}(z)|\leqslant Cs^{2k-1}.
	\end{equation*}
	We therefore scale each eigenvector by $s^{-(2k-1)/2}$ and set
	\begin{equation*}
		u(z)=s^{-(2k-1)/2}\widehat{u}(z),\qquad b(z)=u(z^{-1})^{T}Ju(z).
	\end{equation*}
	Then $\|u(z)\|\leqslant Cs^{-(2k-1)/2}$ and $c\leqslant|b(z)|\leqslant C$. Symplecticity gives
	$u(z)^{T}Ju(z')=0$ unless $zz'=1$, where the pairing is complex bilinear.
	
	For the unit-circle pair, take $\operatorname{Im}z>0$ and write $u(z)=x+iy$. Then
	$b(z)=2it$, where $t=x^{T}Jy$, so
	\begin{equation*}
		\frac{x}{\sqrt{|t|}},\qquad -\frac{\operatorname{sgn}(t)y}{\sqrt{|t|}}
	\end{equation*}
	have pairing $-1$ and give a rotation with $|\psi_{E}|=2\arcsin(s/2)$.
	
	For a stable root $z$, put $v(z)=u(z^{-1})/b(z)$. Then $u(z)^{T}Jv(z)=-1$. Real roots give real
	symplectic pairs. If $z$ is nonreal, symplecticity also gives $u(z)^{T}J\overline{v(z)}=0$.
	Taking real and imaginary parts of these two identities yields
	\begin{equation*}
		\begin{pmatrix}\operatorname{Re}u(z)&\operatorname{Im}u(z)\end{pmatrix}^{T}
		J\begin{pmatrix}\operatorname{Re}v(z)&\operatorname{Im}v(z)\end{pmatrix}
		=\begin{pmatrix}-\frac{1}{2}&0\\0&\frac{1}{2}\end{pmatrix}.
	\end{equation*}
	Thus, for each conjugate pair of nonreal stable roots, we take
	\begin{equation*}
		(\operatorname{Re}u(z),\operatorname{Im}u(z)),\qquad
		(2\operatorname{Re}v(z),-2\operatorname{Im}v(z))
	\end{equation*}
	as the stable and unstable columns, whose cross-pairing matrix is $-I_{2}$. The stable columns
	pair to zero with one another, as do the unstable columns, and pairings between distinct blocks
	vanish. On the stable columns, the matrix is
	\begin{equation*}
		\begin{pmatrix}
			\operatorname{Re}z&\operatorname{Im}z\\
			-\operatorname{Im}z&\operatorname{Re}z
		\end{pmatrix},
	\end{equation*}
	whose norm equals $|z|$. Thus the stable block satisfies $\|D_{E}\|\leqslant e^{-cs}$, and
	symplecticity gives the unstable block $D_{E}^{-T}$. All blocks are $I+O(s)$.
	
	Form $W_{E}$ by placing all stable columns before the corresponding unstable columns. If the
	elliptic pair is present, insert its first vector before the stable columns and its second
	vector before the unstable columns. Then $W_{E}^{T}JW_{E}=J$, with the direct-sum convention
	fixed above. The additional normalization factors are bounded, so each column has norm at most
	$Cs^{-(2k-1)/2}$. This proves the bound for $W_{E}$; the same bound holds for its inverse since
	$W_{E}^{-1}=-JW_{E}^{T}J$.
\end{proof}

\subsection{Block diagonalization}

Let $D\in\operatorname{GL}(d,\mathbb{R})$ and
\begin{equation*}
  H=\operatorname{diag}(D,D^{-T}),\qquad
  \|D\|\leqslant e^{-cs},\qquad \|H-I\|\leqslant C_{0}s,
\end{equation*}
where $c,C_{0}>0$ are fixed and $s>0$ is sufficiently small. The constants in this subsection depend
only on $c,C_{0}$ and the matrix dimension. All matrix functions in this subsection are continuous
and real on $\mathbb{T}$. Rotation numbers are taken with the local real lifts of
Lemma~\ref{l:2.2}(1). An empty block has rotation number zero.

\begin{lemma}\label{l:short-hyperbolic}
If $K$ is symplectic and $\|K-H\|\leqslant c_{1}s$, with $c_{1}>0$ sufficiently small, then
$\rho(\alpha,K)=0$.
\end{lemma}

\begin{proof}
Set $f=K-H$ and $\varepsilon=\|f\|$, and write $f=(f_{ij})_{1\leqslant i,j\leqslant2}$ in the block
decomposition of $H$. We seek an invariant graph $\{(v,X(\theta)v):v\in\mathbb{R}^{d}\}$. Its
invariance is equivalent to
\begin{equation*}
  X(\theta+\alpha)D-D^{-T}X
  =f_{21}+f_{22}X-X(\theta+\alpha)f_{11}-X(\theta+\alpha)f_{12}X.
\end{equation*}
The inverse of the linear operator on the left is
\begin{equation*}
  G\longmapsto-\sum_{n=0}^{\infty}(D^{T})^{n+1}G(\theta+n\alpha)D^{n}.
\end{equation*}
This uniformly convergent series defines an operator of norm at most
$C/(1-e^{-2cs})\leqslant C_{1}/s$. On a ball of radius $r$, the nonlinear right-hand side has norm
at most $\varepsilon(1+r)^{2}$ and Lipschitz constant at most $\varepsilon(2+2r)$. For sufficiently
small $c_{1}$, composition with the linear inverse is therefore a contraction of the closed ball of
radius $r=4C_{1}\varepsilon/s\leqslant1/2$ into itself, with Lipschitz constant at most $1/2$. Its
fixed point satisfies $\|X\|\leqslant C\varepsilon/s$. The induced map on the first coordinate has
norm
\begin{equation*}
  \|D+f_{11}+f_{12}X\|
  \leqslant e^{-cs}+C\varepsilon(1+\varepsilon/s)
  \leqslant e^{-cs/2},
\end{equation*}
after decreasing $c_{1}$. The boundedness of $X$ then gives exponential contraction on the graph.
For any two vectors $v,w$ in the graph, symplecticity gives
\begin{equation*}
  v^{T}J_{2d}w=(K^{(n)}v)^{T}J_{2d}(K^{(n)}w)\longrightarrow0.
\end{equation*}
Thus the graph is Lagrangian and $X=X^{T}$. The symplectic shear
$\begin{pmatrix}I&0\\X(\theta)&I\end{pmatrix}$ is homotopic to the identity and maps the fixed
horizontal Lagrangian subspace onto the graph. Thus Lemma~\ref{l:2.1} and conjugacy invariance give
zero rotation number modulo $\mathbb{Z}$. Since $K=I+O(s)$, its local real lift is zero by
Lemma~\ref{l:2.2}(1).
\end{proof}

The next lemma gives a symplectic block diagonalization with estimates in terms of the hyperbolic
gap $s$; see also \cite{WXZ25}.

\begin{lemma}[Block diagonalization]\label{l:short-center}
Let $M=R_{\psi}\oplus H\in\operatorname{Sp}(2k,\mathbb{R})$, with $|\psi|\leqslant C_{0}s$, and let
$\mathcal{K}=M+f$ be symplectic. If $\varepsilon=\|f\|\leqslant c_{2}s$, with $c_{2}>0$ sufficiently
small, there is a symplectic $Q$, homotopic to the identity, such that
\begin{equation}\label{e:short-center-reduction}
  Q(\theta+\alpha)^{-1}\mathcal{K}(\theta)Q(\theta)
  =\mathcal{C}(\theta)\oplus\widetilde{H}(\theta),\qquad
  \mathcal{C}(\theta)\in\operatorname{SL}(2,\mathbb{R}),
\end{equation}
with
\begin{equation}\label{e:short-center-bounds}
  \|Q-I\|\leqslant C\varepsilon/s,\qquad
  \|\mathcal{C}-R_{\psi}\|+\|\widetilde{H}-H\|\leqslant C\varepsilon.
\end{equation}
Moreover, $\rho(\alpha,\widetilde{H})=0$.
\end{lemma}

\begin{proof}
If the hyperbolic block is absent, take $Q=I$ and $\mathcal{C}=\mathcal{K}$. Otherwise, we first
construct an invariant center graph $\{(v,Z(\theta)v):v\in\mathbb{R}^{2}\}$ and then choose
symplectic bases of this graph and its symplectic orthogonal complement. The graph equation will
give both $\|Z\|=O(\varepsilon/s)$ and $\|Z(\cdot+\alpha)-Z\|=O(\varepsilon)$. The second bound,
together with $M-I=O(s)$, will ensure that the transformed error remains $O(\varepsilon)$.

We work in the block coordinates of the direct sum, with symplectic form
$\widehat{J}=\operatorname{diag}(J_{2},J_{2k-2})$. All matrices in this proof are expressed in these
coordinates. Returning to standard coordinates uses the fixed orthogonal permutation and preserves
the estimates and the homotopy to the identity. Write
\begin{equation*}
  \mathcal{K}(\theta)=\begin{pmatrix}R_{\psi}+f_{11}(\theta)&f_{12}(\theta)\\f_{21}(\theta)&H+f_{22}(\theta)\end{pmatrix}.
\end{equation*}
The four perturbation blocks have norm at most $\varepsilon$. For
$P(\theta)=\binom{I_{2}}{Z(\theta)}$, we have
\begin{equation*}
  \mathcal{K}(\theta)P(\theta)
  =\binom{R_{\psi}+f_{11}+f_{12}Z}{f_{21}+(H+f_{22})Z}.
\end{equation*}
The graph spanned by $P$ is invariant if and only if
$f_{21}+(H+f_{22})Z=Z(\theta+\alpha)(R_{\psi}+f_{11}+f_{12}Z)$, or equivalently
\begin{equation}\label{e:short-center-graph}
  Z(\theta+\alpha)R_{\psi}-HZ
  =f_{21}+f_{22}Z-Z(\theta+\alpha)f_{11}-Z(\theta+\alpha)f_{12}Z.
\end{equation}

Define $\mathcal{L}Z=Z(\cdot+\alpha)R_{\psi}-HZ$ on $C^{0}(\mathbb{T},\mathbb{R}^{(2k-2)\times2})$.
The stable and unstable row blocks of $\mathcal{L}Z=G$ have the solutions
\begin{equation*}
  \begin{aligned}
    Z_{s}(\theta)&=\sum_{n=0}^{\infty}D^{n}G_{s}(\theta-(n+1)\alpha)R_{\psi}^{-(n+1)},\\
    Z_{u}(\theta)&=-\sum_{n=0}^{\infty}(D^{T})^{n+1}G_{u}(\theta+n\alpha)R_{\psi}^{n}.
  \end{aligned}
\end{equation*}
The two series are obtained by backward and forward iteration, respectively. Since $D$ is
contracting and $R_{\psi}$ is orthogonal, they converge uniformly and give the unique bounded
solution, with
\begin{equation*}
  \|\mathcal{L}^{-1}G\|
  \leqslant\frac{C}{1-\|D\|}\|G\|
  \leqslant\frac{C_{1}}{s}\|G\|,
\end{equation*}
where $C_{1}\geqslant1$ is independent of $s$.

Equation~\eqref{e:short-center-graph} is $Z=\Phi(Z)$, where
\begin{equation*}
  \Phi(Z)=\mathcal{L}^{-1}\bigl(f_{21}+f_{22}Z-Z(\cdot+\alpha)f_{11}-Z(\cdot+\alpha)f_{12}Z\bigr).
\end{equation*}
Decrease $c_{2}$ so that $r=4C_{1}\varepsilon/s\leqslant1/2$. For $\|Z\|,\|W\|\leqslant r$,
estimating the nonlinear terms and their differences gives
\begin{equation*}
  \begin{aligned}
    \|\Phi(Z)\|&\leqslant\frac{C_{1}\varepsilon}{s}(1+r)^{2}\leqslant r,\\
    \|\Phi(Z)-\Phi(W)\|&\leqslant\frac{C_{1}\varepsilon}{s}(2+2r)\|Z-W\|
    \leqslant\frac{1}{2}\|Z-W\|.
  \end{aligned}
\end{equation*}
Let $Z$ be the unique fixed point in this ball. The right-hand side of \eqref{e:short-center-graph}
is $O(\varepsilon)$, so
\begin{equation*}
  Z(\theta+\alpha)-Z
  =(H-I)Z-Z(\theta+\alpha)(R_{\psi}-I)+O(\varepsilon).
\end{equation*}
Since $M-I=O(s)$, we obtain
\begin{equation}\label{e:short-graph-variation}
  \|Z\|\leqslant C\varepsilon/s,\qquad
  \|Z(\cdot+\alpha)-Z\|\leqslant C\varepsilon.
\end{equation}
The second estimate will control $Q(\cdot+\alpha)-Q$ after symplectic normalization.

The induced symplectic form on the graph is
\begin{equation*}
  P(\theta)^{T}\widehat{J}P(\theta)
  =J_{2}+Z(\theta)^{T}J_{2k-2}Z(\theta)
  =\omega(\theta)J_{2},\qquad \omega=1+O(\|Z\|^{2}).
\end{equation*}
For sufficiently small $c_{2}$, we have $\omega>0$, so the center graph is symplectic. Its
symplectic orthogonal complement has the frame
\begin{equation*}
  S(\theta)=\binom{J_{2}Z(\theta)^{T}J_{2k-2}}{I_{2k-2}},\qquad
  P(\theta)^{T}\widehat{J}S(\theta)=0.
\end{equation*}
For a center vector $v$ and a vector $w$ in this complement,
\begin{equation*}
  (\mathcal{K}(\theta)v)^{T}\widehat{J}\mathcal{K}(\theta)w=v^{T}\widehat{J}w=0.
\end{equation*}
Since $\mathcal{K}(\theta)$ maps the center graph onto the graph at $\theta+\alpha$, the complement
is invariant.

Divide both columns of $P$ by $\sqrt{\omega}$ and apply symplectic Gram--Schmidt to the columns of
$S$. Pair columns $j$ and $k-1+j$ in increasing order of $j$, preserving their original positions.
At each step, remove the components along all previously normalized pairs by
\begin{equation*}
  x\longmapsto x-(q^{T}\widehat{J}x)p+(p^{T}\widehat{J}x)q,\qquad p^{T}\widehat{J}q=-1,
\end{equation*}
and normalize the resulting pair $(u,v)$ by replacing $v$ with $-v/(u^{T}\widehat{J}v)$. At $Z=0$,
the columns of $S$ are standard coordinate vectors, so each denominator $u^{T}\widehat{J}v$ equals
$-1$; also $\omega=1$. Hence, in a fixed small neighborhood of $Z=0$, independent of $s$, the
denominators stay bounded away from zero and $\omega$ stays positive. Decrease $c_{2}$ so that
$Z(\theta)$ stays in this neighborhood. The resulting completion
$\mathcal{Q}(Z)$ is therefore smooth, with $\mathcal{Q}(0)=I$ and uniformly bounded derivative. The columns of
$Q(\theta)=\mathcal{Q}(Z(\theta))$ form symplectic bases of the two invariant subspaces, which gives
\eqref{e:short-center-reduction}. The center block is
\begin{equation*}
  \mathcal{C}(\theta)=\sqrt{\frac{\omega(\theta+\alpha)}{\omega(\theta)}}
  \bigl(R_{\psi}+f_{11}(\theta)+f_{12}(\theta)Z(\theta)\bigr)
  \in\operatorname{SL}(2,\mathbb{R}).
\end{equation*}
The bounded derivative of $\mathcal{Q}$ and \eqref{e:short-graph-variation} give
\begin{equation*}
  \|Q-I\|\leqslant C\varepsilon/s,\qquad
  \|Q(\cdot+\alpha)-Q\|\leqslant C\varepsilon.
\end{equation*}
The path $\mathcal{Q}(tZ(\theta))$, $0\leqslant t\leqslant1$, is a homotopy from the identity to
$Q$.

The diagonal blocks are estimated from
\begin{equation*}
  \begin{aligned}
    Q(\theta+\alpha)^{-1}(M+f)Q-M
    =Q(\theta+\alpha)^{-1}\bigl(&[M,Q]\\
    &-(Q(\theta+\alpha)-Q)M+fQ\bigr),
  \end{aligned}
\end{equation*}
where $[M,Q]=MQ-QM$. Since
\begin{equation*}
  \|[M,Q]\|=\|[M-I,Q-I]\|
  \leqslant2\|M-I\|\|Q-I\|\leqslant C\varepsilon,
\end{equation*}
each term in brackets is $O(\varepsilon)$. The uniform bound on $Q^{-1}$ then gives
\eqref{e:short-center-bounds}. For small $c_{2}$, Lemma~\ref{l:short-hyperbolic} applies to
$\widetilde{H}$ and gives $\rho(\alpha,\widetilde{H})=0$.
\end{proof}

\subsection{Rotation number and IDS estimates}

Fix $p_{*}$ and $B_{*}$ from Proposition~\ref{p:4.1}, with $h_{0}=1$ and $h_{*}=1/2$, and write
$A_{E}(\theta)=A_{E}(p_{*},\theta)$. All constants below depend only on $k,\gamma,\tau$ and are
uniform for $\lambda\geqslant\lambda_{0}$. We may increase $\lambda_{0}$ and decrease the radius
$E_{*}$ in Lemma~\ref{l:short-reference}.

The coordinates of Section~2 give
\begin{equation*}
  \partial_{E}A_{E}(p,\theta)=e_{1}e_{k}^{T},
\end{equation*}
because $\partial_{E}S_{E}=e_{1}e_{k}^{T}$, $\mathcal{T}(p)^{-1}e_{1}=e_{1}$, and
$e_{k}^{T}\mathcal{T}(p)=e_{k}^{T}$. Together with \eqref{e:4.1}, this gives
\begin{equation*}
  \mathcal{A}_{E}(\theta)
  :=B_{*}(\theta+\alpha)^{-1}A_{E}(\theta)B_{*}(\theta)
  =U_{*}+E B_{*}(\theta+\alpha)^{-1}e_{1}e_{k}^{T}B_{*}(\theta).
\end{equation*}
Since $\bar{A}_{E}(\bar{p})=U_{*}+Ee_{1}e_{k}^{T}$, estimate~\eqref{e:4.2} implies
\begin{equation}\label{e:short-reference-error}
  \|\mathcal{A}_{E}-\bar{A}_{E}(\bar{p})\|
  \leqslant C\lambda^{-1}|E|.
\end{equation}

For $0<|E|<E_{*}$, put $s=|E|^{1/(2k)}$ and use $W_{E}$ from Lemma~\ref{l:short-reference}. By
\eqref{e:short-reference-error} and \eqref{e:short-condition},
\begin{equation}\label{e:short-scaled-error}
  W_{E}^{-1}\mathcal{A}_{E}(\theta)W_{E}
  =\bar{M}_{E}+f_{E}(\theta),\qquad
  \|f_{E}\|\leqslant C\lambda^{-1}s^{2k}s^{-(2k-1)}=C\lambda^{-1}s.
\end{equation}

The matrix $U_{*}=e^{N}$ preserves the Lagrangian subspace $\operatorname{Im}N^{k}$, since this
subspace has dimension $k$ and
\begin{equation*}
  (N^{k})^{T}JN^{k}=(-1)^{k}JN^{2k}=0.
\end{equation*}
Thus $\rho(U_{*})\in\mathbb{Z}$ by Lemma~\ref{l:2.1}. Since $B_{*}$ is homotopic to the identity, we
may choose the continuous real lift with $\rho(\alpha,A_{0})=0$.

Equations \eqref{e:short-model-gap} and \eqref{e:short-scaled-error} imply
\begin{equation*}
  \|W_{E}^{-1}\mathcal{A}_{E}W_{E}-I\|\leqslant Cs.
\end{equation*}
Lemma~\ref{l:2.2}(1) and conjugacy invariance therefore give
\begin{equation*}
  \operatorname{dist}\bigl(\rho(\alpha,A_{E}),\mathbb{Z}\bigr)
  \leqslant C|E|^{1/(2k)}.
\end{equation*}
Conjugacy by $W_{E}$ is used here only at each fixed $E\ne0$. The chosen real lift is continuous
along the original family $E\mapsto A_{E}$ and equals zero at $E=0$. After decreasing $E_{*}$ so that
the right-hand side is less than $1/8$, continuity gives
\begin{equation}\label{e:short-uniform-lift}
  |\rho(\alpha,A_{E})|<1/8,\qquad |E|<E_{*},
\end{equation}
uniformly in $\lambda$ and $\alpha$ under consideration.

\begin{proposition}\label{p:short-rotation}
For $0<|E|<E_{*}$ with $(-1)^{k}E>0$,
\begin{equation}\label{e:short-interior-rho}
  \left|\,|\rho(\alpha,A_{E})|
  -\frac{2}{\pi}\arcsin\!\left(\frac{|E|^{1/(2k)}}{2}\right)\right|
  \leqslant C\lambda^{-1}|E|^{1/(2k)}.
\end{equation}
For $0<|E|<E_{*}$ with $(-1)^{k}E<0$,
\begin{equation}\label{e:short-exterior-rho}
  \rho(\alpha,A_{E})=0.
\end{equation}
\end{proposition}

\begin{proof}
For sufficiently large $\lambda_{0}$, Lemmas~\ref{l:short-hyperbolic} and~\ref{l:short-center} apply
to \eqref{e:short-scaled-error}. If $(-1)^{k}E>0$, Lemma~\ref{l:short-center} gives the direct sum
$\mathcal{C}_{E}\oplus\widetilde{H}_{E}$, with
\begin{equation*}
  \|\mathcal{C}_{E}-R_{\psi_{E}}\|\leqslant C\lambda^{-1}s,\qquad
  \rho(\alpha,\widetilde{H}_{E})=0.
\end{equation*}
Lemma~\ref{l:2.2}(2) yields
\begin{equation*}
  \rho(\alpha,\mathcal{C}_{E})=\frac{\psi_{E}}{\pi}+O(\lambda^{-1}s).
\end{equation*}
By additivity and conjugacy invariance, this equals $\rho(\alpha,A_{E})$ modulo $\mathbb{Z}$. For
sufficiently small $E_{*}$, both real lifts have absolute value less than $1/8$, so they agree.
Taking absolute values and using $|\psi_{E}|=2\arcsin(s/2)$ proves \eqref{e:short-interior-rho}.

If $(-1)^{k}E<0$, Lemma~\ref{l:short-hyperbolic} applies to the full cocycle in
\eqref{e:short-scaled-error}. The rotation number is zero modulo $\mathbb{Z}$, and hence zero for
the chosen lift by \eqref{e:short-uniform-lift}.
\end{proof}

\begin{proof}[Proof of Theorem~\ref{t:dual}]
By \eqref{e:short-exterior-rho} and \eqref{e:2.13}, the IDS is constant on $(-1)^{k+1}(0,E_{*})$.
Continuity gives the same value at zero, and \eqref{e:2.14} excludes this interval from the
spectrum. Let $0<r<E_{*}$. If $k$ is even, then
$\widetilde{N}_{p_{*}}(-r)=\widetilde{N}_{p_{*}}(0)$, and \eqref{e:2.13} gives
\begin{equation*}
  \widetilde{N}_{p_{*}}(r)-\widetilde{N}_{p_{*}}(-r)
  =-\rho(\alpha,A_{r}).
\end{equation*}
If $k$ is odd, then $\widetilde{N}_{p_{*}}(r)=\widetilde{N}_{p_{*}}(0)$, and the same identity gives
\begin{equation*}
  \widetilde{N}_{p_{*}}(r)-\widetilde{N}_{p_{*}}(-r)
  =\rho(\alpha,A_{-r}).
\end{equation*}
Since the IDS is nondecreasing, both cases yield
\begin{equation*}
  \widetilde{N}_{p_{*}}(r)-\widetilde{N}_{p_{*}}(-r)
  =|\rho(\alpha,A_{(-1)^{k}r})|.
\end{equation*}
Proposition~\ref{p:short-rotation} therefore yields
\begin{equation*}
  \left|\widetilde{N}_{p_{*}}(r)-\widetilde{N}_{p_{*}}(-r)
  -\frac{2}{\pi}\arcsin\!\left(\frac{r^{1/(2k)}}{2}\right)\right|
  \leqslant C\lambda^{-1}r^{1/(2k)}.
\end{equation*}
Since $2\arcsin(s/2)=s+O(s^{3})$, we obtain, after increasing $\lambda_{0}$ and decreasing $E_{*}$,
\begin{equation*}
  cr^{1/(2k)}\leqslant\widetilde{N}_{p_{*}}(r)-\widetilde{N}_{p_{*}}(-r)
  \leqslant Cr^{1/(2k)},\qquad 0<r<E_{*}.
\end{equation*}
This shows that $0\in\widetilde{\Sigma}_{p_{*}}$ by \eqref{e:2.14}.

It remains to show that zero is an outer spectral edge. By \eqref{e:4.2}, uniformly in $\theta$,
\begin{equation*}
  \|L_{p_{*},\theta}-\bar{L}_{\bar{p}}\|
  \leqslant\sum_{j=1}^{k-1}|a_{j,*}-\bar{a}_{j}|+|e_{*}-\bar{e}|+2\lambda^{-1}
  \leqslant C\lambda^{-1}.
\end{equation*}
The reference spectrum is $\operatorname{spec}(\bar{L}_{\bar{p}})=(-1)^{k}[0,4^{k}]$. Hence any
spectrum on the exterior side lies in $(-1)^{k+1}(0,C\lambda^{-1}]$. With $E_{*}$ fixed, increase
$\lambda_{0}$ so that $C\lambda_{0}^{-1}<E_{*}$. This interval is then contained in the excluded
interval. Thus zero is the bottom of $\widetilde{\Sigma}_{p_{*}}$ for even $k$ and the top for odd
$k$. The parameter bound follows from \eqref{e:4.2}.
\end{proof}

\begin{proof}[Proof of Theorem~\ref{t:main}]
	Take $\lambda>\max\{1,\lambda_{0}\}$ and let $p_{*}$ be the parameter constructed in the proof of
	Theorem~\ref{t:dual}. Set $v=\lambda V_{p_{*}}$ and $E_{0}=0$. By \eqref{e:au}, zero is an outer
	edge of $\Sigma_{v}$. The linear rescaling of energy preserves the local exponents, which equal
	$1/(2k)$ by Theorem~\ref{t:dual}. Corollary~\ref{c:4.4} gives $L_{0}^{v}(0)=\log\lambda>0$ and
	$\omega^{v}(0)=k$.
\end{proof}

\section*{Acknowledgements}
Jiangong You is supported by NSFC grants (12531006, 12526201) and Nankai Zhide Foundation.

\section*{Declaration of AI use}
The authors used ChatGPT to improve the language and to help identify potential errors in the
mathematical arguments. All suggestions were independently assessed by the authors, who take full
responsibility for the manuscript.

\end{document}